\documentclass[10pt, reqno]{amsart}
\pdfoutput=1
\usepackage{appendix}
\usepackage{amsmath}
\usepackage[a4paper, centering]{geometry}
\usepackage{xcolor}
\usepackage[utf8]{inputenc}
\usepackage[colorlinks=true,hyperindex,pagebackref,linktocpage=true]{hyperref}
\usepackage{cleveref}
\usepackage{mathtools}
\usepackage{amsfonts}
\usepackage{amssymb}
\usepackage{tikz-cd}
\usepackage[T1]{fontenc}

\hypersetup{
	colorlinks,
	linkcolor={violet!60!black},
	citecolor={cyan!60!black},
	urlcolor={orange!60!black}
}

\usepackage{stmaryrd}
\usepackage{mathrsfs}  
\usepackage{varwidth}
\theoremstyle{plain}
\newtheorem{theorem}{Theorem}[section]
\newtheorem{proposition}[theorem]{Proposition}
\newtheorem{lemma}[theorem]{Lemma}
\newtheorem{corollary}[theorem]{Corollary}
\newtheorem{conjecture}[theorem]{Conjecture}
\newtheorem{question}[theorem]{Question}

\theoremstyle{definition}
\newtheorem{definition}[theorem]{Definition}

\newtheorem{remark}[theorem]{Remark}

\newcommand{\nc}{\newcommand}
\nc{\on}{\operatorname}

\nc{\Q}{\mathbb{Q}}
\nc{\Z}{\mathbb{Z}}
\nc{\cl}{\mathrm{cl}}

\nc{\fraka}{{\mathfrak a}} \nc{\bba}{{\mathbf a}}
\nc{\frakb}{{\mathfrak b}}
\nc{\frakc}{{\mathfrak c}}
\nc{\frakd}{{\mathfrak d}}
\nc{\frake}{{\mathfrak e}}
\nc{\frakf}{{\mathfrak f}}
\nc{\frakg}{{\mathfrak g}}
\nc{\frakh}{{\mathfrak h}}
\nc{\fraki}{{\mathfrak i}}
\nc{\frakj}{{\mathfrak j}}
\nc{\frakk}{{\mathfrak k}}
\nc{\frakl}{{\mathfrak l}}
\nc{\frakm}{{\mathfrak m}}
\nc{\frakn}{{\mathfrak n}}
\nc{\frako}{{\mathfrak o}}
\nc{\frakp}{{\mathfrak p}}
\nc{\frakq}{{\mathfrak q}}
\nc{\frakr}{{\mathfrak r}}
\nc{\fraks}{{\mathfrak s}}
\nc{\frakt}{{\mathfrak t}}
\nc{\fraku}{{\mathfrak u}}
\nc{\frakv}{{\mathfrak v}}
\nc{\frakw}{{\mathfrak w}}
\nc{\frakx}{{\mathfrak x}}
\nc{\fraky}{{\mathfrak y}}
\nc{\frakz}{{\mathfrak z}}
\nc{\frakA}{{\mathfrak A}}
\nc{\frakB}{{\mathfrak B}}
\nc{\frakC}{{\mathfrak C}}
\nc{\frakD}{{\mathfrak D}}
\nc{\frakE}{{\mathfrak E}}
\nc{\frakF}{{\mathfrak F}}
\nc{\frakG}{{\mathfrak G}}
\nc{\frakH}{{\mathfrak H}}
\nc{\frakI}{{\mathfrak I}}
\nc{\frakJ}{{\mathfrak J}}
\nc{\frakK}{{\mathfrak K}}
\nc{\frakL}{{\mathfrak L}}
\nc{\frakM}{{\mathfrak M}}
\nc{\frakN}{{\mathfrak N}}
\nc{\frakO}{{\mathfrak O}}
\nc{\frakP}{{\mathfrak P}}
\nc{\frakQ}{{\mathfrak Q}}
\nc{\frakR}{{\mathfrak R}}
\nc{\frakS}{{\mathfrak S}}
\nc{\frakT}{{\mathfrak T}}
\nc{\frakU}{{\mathfrak U}}
\nc{\frakV}{{\mathfrak V}}
\nc{\frakW}{{\mathfrak W}}
\nc{\frakX}{{\mathfrak X}}
\nc{\frakY}{{\mathfrak Y}}
\nc{\frakZ}{{\mathfrak Z}}
\nc{\bbA}{{\mathbb A}}
\nc{\bbC}{{\mathbb C}}
\nc{\bbD}{{\mathbb D}}
\nc{\bbE}{{\mathbb E}}
\nc{\bbF}{{\mathbb F}} \nc{\bbf}{{\mathbf f}}
\nc{\bbG}{{\mathbb G}}
\nc{\bbH}{{\mathbb H}}
\nc{\bbI}{{\mathbb I}}
\nc{\bbJ}{{\mathbb J}}
\nc{\bbK}{{\mathbb K}}
\nc{\bbL}{{\mathbb L}}
\nc{\bbM}{{\mathbb M}}
\nc{\bbN}{{\mathbb N}}
\nc{\bbO}{{\mathbb O}}
\nc{\bbP}{{\mathbb P}}
\nc{\bbQ}{{\mathbb Q}}
\nc{\bbR}{{\mathbb R}}
\nc{\bbS}{{\mathbb S}}
\nc{\bbT}{{\mathbb T}}
\nc{\bbU}{{\mathbb U}}
\nc{\bbV}{{\mathbb V}}
\nc{\bbW}{{\mathbb W}}
\nc{\bbX}{{\mathbb X}}
\nc{\bbY}{{\mathbb Y}}
\nc{\bbZ}{{\mathbb Z}}
\nc{\calA}{{\mathcal A}}
\nc{\calB}{{\mathcal B}}
\nc{\calC}{{\mathcal C}}
\nc{\calD}{{\mathcal D}}
\nc{\calE}{{\mathcal E}}
\nc{\calF}{{\mathcal F}}
\nc{\calG}{{\mathcal G}}
\nc{\calH}{{\mathcal H}}
\nc{\calI}{{\mathcal I}}
\nc{\calJ}{{\mathcal J}}
\nc{\calK}{{\mathcal K}}
\nc{\calL}{{\mathcal L}}
\nc{\calM}{{\mathcal M}}
\nc{\calN}{{\mathcal N}}
\nc{\calO}{{\mathcal O}}
\nc{\calP}{{\mathcal P}}
\nc{\calQ}{{\mathcal Q}}
\nc{\calR}{{\mathcal R}}
\nc{\calS}{{\mathcal S}}
\nc{\calT}{{\mathcal T}}
\nc{\calU}{{\mathcal U}}
\nc{\calV}{{\mathcal V}}
\nc{\calW}{{\mathcal W}}
\nc{\calX}{{\mathcal X}}
\nc{\calY}{{\mathcal Y}}
\nc{\calZ}{{\mathcal Z}}

\nc{\scrA}{{\mathscr A}}
\nc{\scrB}{{\mathscr B}}
\nc{\scrC}{{\mathscr C}}
\nc{\scrD}{{\mathscr D}}
\nc{\scrE}{{\mathscr E}}
\nc{\scrF}{{\mathscr F}}
\nc{\scrG}{{\mathscr G}}
\nc{\scrH}{{\mathscr H}}
\nc{\scrI}{{\mathscr I}}
\nc{\scrJ}{{\mathscr J}}
\nc{\scrK}{{\mathscr K}}
\nc{\scrL}{{\mathscr L}}
\nc{\scrM}{{\mathscr M}}
\nc{\scrN}{{\mathscr N}}
\nc{\scrO}{{\mathscr O}}
\nc{\scrP}{{\mathscr P}}
\nc{\scrQ}{{\mathscr Q}}
\nc{\scrR}{{\mathscr R}}
\nc{\D}{{\on{D}}}
\nc{\Div}{{\on{Div}}}
\nc{\Perv}{{\on{Perv}}}

\nc{\bnu}{{\bar{ \nu}}}
\nc{\olO}{\bar{\calO}}

\nc{\sh}{\mathrm{sh}} 
\nc{\al}{{\alpha}} 
\nc{\be}{{\beta}}
\nc{\ga}{{\gamma}} \nc{\Ga}{{\Gamma}}
\nc{\hGa}{\hat{\Gamma}}
\nc{\ve}{{\varepsilon}} 
\nc{\la}{{\lambda}} \nc{\La}{{\Lambda}}
\nc{\om}{\omega} \nc{\Om}{\Omega} 
\nc{\sig}{{\sigma}} \nc{\Sig}{{\Sigma}}
\nc{\dR}{{\mathrm{dR}}}
\nc{\Perf}{{\mathrm{Perf}}}
\nc{\perf}{{\mathrm{perf}}}
\nc{\Gm}{{\mathbb{G}_m}}
\nc{\colim}{{\on{colim}}}
\nc{\et}{\mathrm{\acute{e}t}}
\nc{\Proj}{{\mathrm{Proj}}}

\makeatletter
\DeclareFontEncoding{LS1}{}{}
\DeclareFontSubstitution{LS1}{stix}{m}{n}
\DeclareMathAlphabet{\rhomalpha}{LS1}{stixscr}{m}{n}
\makeatother

\nc{\Spa}{\on{{Spa}}}
\nc{\Spd}{\on{{Spd}}}
\nc{\tnb}{\psi_{\rm tame}}
\nc{\oM}{\overline{{M}}}
\nc{\op}{{\on{op}}}
\nc{\ad}{{\on{ad}}}
\nc{\alg}{{\on{alg}}}
\nc{\Ad}{{\on{Ad}}}
\nc{\Adm}{{\on{Adm}}} \nc{\aff}{{\on{af}}}
\nc{\Aut}{{\on{Aut}}}
\nc{\Bun}{{\on{Bun}}}
\nc{\cha}{{\on{char}}}
\nc{\der}{{\on{der}}}
\nc{\Der}{{\on{Der}}}
\nc{\diag}{{\on{diag}}}
\nc{\End}{{\on{End}}}
\nc{\Fl}{{\mathrm{Fl}}}
\nc{\Tr}{{\on{Transp}}}
\nc{\TR}{{\calT\!\calR}}
\nc{\Gal}{{\on{Gal}}}
\nc{\Gr}{{\on{Gr}}}
\nc{\Hk}{{\on{Hk}}}
\nc{\rH}{{\on{H}}}
\nc{\Hom}{{\on{Hom}}}
\nc{\IC}{{\on{IC}}}
\nc{\id}{{\on{id}}}
\nc{\Id}{{\on{Id}}}
\nc{\ind}{{\on{ind}}}
\nc{\Ind}{{\on{Ind}}}
\nc{\Lie}{{\on{Lie}}}
\nc{\Pic}{{\on{Pic}}}
\nc{\pr}{{\on{pr}}}
\nc{\Res}{{\on{Res}}}
\nc{\res}{{\on{res}}} \nc{\Sat}{{\on{Sat}}}
\nc{\spc}{{\on{sc}}}
\nc{\drv}{{\on{der}}}
\nc{\sgn}{{\on{sgn}}}
\nc{\Spec}{{\on{Spec}}}\nc{\Spf}{\on{Spf}} 
\nc{\Sph}{\on{Sph}}
\nc{\St}{{\on{St}}}
\nc{\tr}{{\on{tr}}}
\nc{\Mod}{{\mathrm{-Mod}}}
\nc{\Hilb}{{\on{Hilb}}} 
\nc{\Ext}{{\on{Ext}}} 
\nc{\vs}{{\on{Vec}}}
\nc{\ev}{{\on{ev}}}
\nc{\nO}{{\breve{\calO}}}
\nc{\tS}{{\tilde{S}}}
\nc{\spe}{{\on{sp}}}
\nc{\loc}{{\on{loc}}}
\nc{\pre}{{\on{pre}}}
\nc{\Alg}{{\on{Alg}_k^{\on{pf}}}}

\nc{\dimt}{{\on{dim.trg}}}

\nc{\co}{\colon}
\nc{\dia}{{\diamondsuit}}

\nc{\nscrR}{{\mathscr{R}^{\on{nr}}}}

\nc{\GL}{{\on{GL}}}

\nc{\Gl}{\on{Gl}} 
\nc{\GSp}{{\on{GSp}}}
\nc{\gl}{{\frakg\frakl}}
\nc{\SL}{{\on{SL}}} 
\nc{\SU}{{\on{SU}}} 
\nc{\SO}{{\on{SO}}}
\nc{\PGL}{{\on{PGL}}}

\nc{\Conv}{{\on{Conv}}}
\nc{\Rep}{{\on{Rep}}}
\nc{\Dom}{{\on{Dom}}}
\nc{\red}{{\on{red}}}
\nc{\act}{{\on{act}}}
\nc{\nr}{{\on{nr}}}
\nc{\ctf}{{\on{ctf}}}

\nc{\str}{{\on{-}}} 
\nc{\os}{{\bar{s}}}
\nc{\oeta}{{\bar{\eta}}}

\nc{\hookto}{\hookrightarrow}
\nc{\longto}{\longrightarrow}
\nc{\leftto}{\leftarrow}
\nc{\onto}{\twoheadrightarrow}
\nc{\lonto}{\twoheadleftarrow}

\numberwithin{equation}{section}

\begin{document}
	
	\title{Mod $p$ sheaves on Witt flags}
	
    \author[R.~Cass]{Robert Cass}
    \address{Mathematical Sciences Department, Claremont McKenna College, 850 Columbia Avenue, Claremont, CA 91711 USA}
    \email{robert.cass@claremontmckenna.edu}

   \author[J.~Louren\c{c}o]{Jo\~ao Louren\c{c}o}
 \address{L'Institut Galilée, Université Sorbonne Paris Nord, 99 avenue Jean Clément, Villetaneuse, France}
 \email{lourenco@math.univ-paris13.fr}

	\begin{abstract}
     We characterize Cohen--Macaulay and $\varphi$-rational perfect schemes in terms of their perverse étale $\bbF_p$-sheaves.
     Using inversion of adjunction, we prove that sufficiently small Schubert varieties in the Witt affine flag variety are perfections of globally $+$-regular varieties, and hence they are $\varphi$-rational. Our methods apply uniformly to all affine Schubert varieties in equicharacteristic, as well as classical Schubert varieties, thereby answering a question of Bhatt. As a corollary, we deduce that scheme-theoretic local models always have $\varphi$-split special fiber.
	\end{abstract}

	\maketitle
	\tableofcontents
	
\section{Introduction}
Let $F$ be a non-archimedean local field with perfect residue field $k$ of characteristic $p>0$. Categories of sheaves on affine flag varieties of reductive groups over $F$ play a major role in geometric representation theory and in geometric approaches to the Langlands program. Here the relevant geometry happens over the residue field $k$. If one considers $\bbF_\ell$- or $\mathbb{Q}_\ell$-étale sheaves for some prime $\ell\neq p$, then these categories are well-studied for $F$ of equal or mixed characteristic. For example, one has the geometric Satake equivalence, e.g.~\cite{MV07, Zhu17}, a collection of central sheaves, e.g.~\cite{Gai01, ALWY23}, and even $\mathbb{Q}$- and $\mathbb{Z}$-linear motivic sheaves, e.g.~\cite{RS21, CvdHS25, CvdHS24, vdH24}.

The situation changes drastically when one considers étale sheaves for $\ell=p$. There is still a perverse t-structure on $k$-varieties due to Gabber \cite{Gab04}, but its behavior can be quite strange, as half of the six functors do not preserve constructibility. When $F$ has characteristic $p$, the first author studied perverse $\mathbb{F}_p$-sheaves on affine flag varieties in \cite{Cas21, Cas22}. There a decisive role is played by the relationship between singularities and  the Frobenius morphism. In this paper we denote the absolute Frobenius morphism on $k$-schemes by $\varphi$, reserving $F$ for our local field.

The purpose of this paper, as the title suggests, is to launch an investigation of perversity properties when $F$ has characteristic $0$. We divide this introduction into four parts: in the first two we explain our purely algebro-geometric results, and in the last two our applications to varieties occurring in the Langlands program.

\subsection{Mod $p$ intersection cohomology}
When $F$ has characteristic $0$, affine flag varieties exist only canonically as functors on perfect $k$-algebras. In \cite{CX25}, some applications of $\mathbb{F}_p$-sheaves to the mod $p$ Langlands program over $F$ were obtained, but with no concern for the perverse t-structure.
Thus, our first order of business is to investigate IC sheaves. Toward this direction, we prove the following general result.

\begin{theorem} \label{thm1}
	Let $k$ be a perfect field of characteristic $p$ and let $X$ be a connected perfectly finite type $k$-scheme. Then $X$ is Cohen--Macaulay (resp.~$\varphi$-rational) if and only if the shifted constant sheaf $\bbF_p[\dim X]$ is perverse (resp.~perverse and simple).
\end{theorem}

To explain the notions in the above theorem, recall that for a local ring $(R,\mathfrak{m})$ of characteristic $p$, the local cohomology groups $H^i_{\mathfrak{m}}(R)$ are modules over the non-commutative polynomial ring $R[\varphi]$.
We say that a perfect local ring $R$ is Cohen--Macaulay if and only if the local cohomology groups $H^i_{\mathfrak{m}}(R)$ vanish for $i < \dim R$, and $\varphi$-rational if, in addition, $H^{\dim R}_{\mathfrak{m}}(R)$ is a simple $R[\varphi]$-module. 
Passing to a deperfection $R_0$, this recovers the classical notions up to $\varphi$-torsion: we prove in Lemma \ref{lemm--rationalequiv} that $R$ is $\varphi$-rational if and only if $R_0$ is $\varphi$-nilpotent. The latter property is a topic of active research in commutative algebra, see e.g.~\cite{ST17, PQ19, DMP24, KMPS23}.

In \cite{Cas22} it was shown that the geometric properties in Theorem \ref{thm1} for a finite type $k$-scheme imply the corresponding properties of perverse $\mathbb{F}_p$-sheaves, but the fact that the converse holds after passing to the perfection lies much deeper. As it turns out, Theorem \ref{thm1} was known to  experts in the $\varphi$-singularities community and appeared in \cite{BBL+23}, after we already found an argument independently. We have included our argument for the benefit of readers unfamiliar with the literature on  $\varphi$-singularities, and because it differs significantly from the one in \cite{BBL+23} in that we perform most of the key arguments on the coherent as opposed to topological side.

\subsection{Abstract Bott--Samelson--Demazure--Hansen varieties}

The property of $\varphi$-rationality is of a local nature, and in particular, it does not descend along proper covers. In order to get proper descent, one has to define a global variant of $\varphi$-rationality, but it is unclear how to proceed in the perfect setting. For noetherian schemes, this is well understood via the property of strong $\varphi$-regularity of Hochster--Huneke \cite{HH89} and its global variant \cite{Smi00}.

Previous results on the $\varphi$-singularities of Schubert varieties were deduced from 
applying the Mehta--Ramanathan criterion, see \cite{MR85}, on Bott--Samelson--Demazure--Hansen (BSDH) resolutions. This presupposes the existence of certain integral divisors which are delicate to construct, see e.g.~\cite{Fal03}. In this paper, we use the closely related property of global $+$-regularity introduced in \cite{BMP+23}, as it carries the advantage of making every $\bbQ$-divisor integral up to passing to a cyclic cover. 

We define an (abstract) BSDH $k$-variety $X$ to be a successive locally Zariski $\mathbb{P}^1_k$-fibration equipped with sections, starting with the base case of a point. It is possible to identify the Picard group with  $\bbZ^{\oplus \mathrm{dim}X}$ by restriction to the distinguished copies of $\mathbb{P}^1_k$, and also to define the boundary divisor $\partial X$. We say the BSDH $k$-variety $X$ is based if it admits an origin-avoiding effective $\bbQ$-divisor $\theta_X$. A way to produce this is via some global section $\vartheta_{X,n} \in H^0(X,\calO(n,\ldots,n))$ not vanishing at the origin for sufficiently divisible $n$ and take $\theta_X=n^{-1}\mathrm{div}(\vartheta_{X,n})$. Using inversion of adjunction from \cite[Theorem 7.2]{BMP+23} (recalled in simplified form in Theorem \ref{thm_inversion_adjunction}), we prove the following result. 

\begin{theorem} \label{BSDHTheorem}
    Let $X$ be a based BSDH $k$-variety and $\theta_X$ an origin-avoiding effective $\bbQ$-divisor. Then the pair $(X,\partial X+\theta_X)$ is asymptotically globally $+$-regular.
\end{theorem}

Asymptotic global $+$-regularity means that $(X,\Delta)$ is globally $+$-regular in the sense of Definition \ref{def:G+R} for every effective $\bbQ$-divisor $\Delta \leq (1-\epsilon)(\partial_X + \theta_X)$ for all $\epsilon>0$, and it implies $\varphi$-rationality. The idea of the proof is to slightly perturb the coefficients of the boundary divisor in such a way that the anti-canonical divisor of the pair becomes ample. As a corollary, we get asymptotic global $\varphi$-regularity of $X$, and compatible $\varphi$-splitting with BSDH subvarieties.
Additionally, we obtain global $+$-regularity for Stein factorizations of semi-ample line bundles on $X$. There is also a perfect analog of BSDH $k$-varieties for which one can ask whether $\varphi$-rationality holds, but there is less evidence for this than in the group-theoretic setting of the next section.

\subsection{Affine Schubert varieties}
Let $\mathcal{G}$ be a parahoric model of a connected reductive group $G$ over $F$. Recall that the affine flag variety $\Fl_{\mathcal{G}}$ is defined as the étale quotient of the loop group $LG$ by the positive loop group $L^+\calG$. In order to shorten our notation, we denote these loop groups by $\mathsf{G}$ and $\mathsf{P}$, respectively, and the flag variety by $\mathsf{G}/\mathsf{P}$. The Bruhat decomposition yields perfectly projective Schubert schemes $\mathsf{S}_w \subset \mathsf{G}/\mathsf{P}$ indexed by double cosets $w$ of the Iwahori--Weyl group.
When $F$ has characteristic $p$ each $\mathsf{S}_{w}$ is canonically isomorphic to the perfection of a projective $k$-scheme (the seminormalization of the affine Schubert variety in \cite{PR08}, see also \cite{HLR24,FHLR25} for the necessity of this functor). It was proved in \cite{Cas22} for split $G$ and \cite{FHLR25} for almost all $G$ that these deperfections are globally $\varphi$-regular. When $F$ has characteristic $0$ we make the following conjecture, which we prove in some cases.
\begin{conjecture}
	The perfect Schubert schemes $\mathsf{S}_{w}$ are $\varphi$-rational.
\end{conjecture}

To bring BSDH $k$-varieties into the picture,  
recall that the multiplication map of $\mathsf{G}$ affords a convolution structure on the Hecke stack $\mathsf{P}\backslash \mathsf{G}/\mathsf{P}$. Forgetting the left $\mathsf{P}$-quotient yields partial resolutions of $\mathsf{S}_w$. Concretely, for any reduced word $s_{\bullet}$ for $w$ we have the Demazure variety
\begin{equation} \label{Demeq}
\mathsf{S}_{s_\bullet}:=\mathsf{P}_{s_1}\times^{\mathsf{B}}\dots\times^{\mathsf{B}}\mathsf{P}_{s_n}/\mathsf{B}
\end{equation}
for an adequate choice of a pro-solvable subgroup $\mathsf{B}\subset \mathsf{P}$. The equation \eqref{Demeq} exhibits $\mathsf{S}_{s_\bullet}$ as an iterated $\mathbb{P}^{1, \mathrm{pf}}_k$-bundle with canonical sections. 
One can show that $\mathsf{S}_{s_{\bullet}}$ has a canonical deperfection by a based BSDH $k$-variety $\mathsf{S}_{s_{\bullet}}^{\mathrm{can}}$ under the assumption that the $\mathsf{P}$-action on $\mathsf{S}_w$ factors through the mod $p$ fiber of $\calG$.
We denote the boundary divisor by $\partial\mathsf{S}_{s_{\bullet}}^{\mathrm{can}}$ 
and the theta divisor $\theta_{s_\bullet}^{\mathrm{can}}$ can be chosen arbitrarily.
The following result, obtained from applying Theorem \ref{BSDHTheorem}, is as sharp as the ones found in the literature for $\varphi$-splittings, e.g.~compare with Lauritzen--Raben-Pedersen--Thomsen \cite{LRPT06} in the case of finite flag varieties.

\begin{theorem} \label{thm2}
	Assume $s_\bullet$ is a reduced word for $w$ and the $\mathsf{P}$-action on $\mathsf{S}_w$ factors through the mod $p$ fiber of $\calG$. Then $(\mathsf{S}_{s_{\bullet}}^{\mathrm{can}},\partial\mathsf{S}_{s_{\bullet}}^{\mathrm{can}}+\theta_{s_\bullet}^{\mathrm{can}})$ is asymptotically globally $+$-regular.
\end{theorem}

As a corollary, whenever the hypotheses of Theorem \ref{thm2} are satisfied, we deduce that the Stein factorization $\mathsf{S}_{w}^{\textrm{can}}$ is globally $\varphi$-regular and compatibly $\varphi$-split with all Schubert subvarieties. The hypothesis on the $\mathsf{P}$-action is trivially satisfied when $F$ has characteristic $p$, in which case we obtain a new proof of the global $\varphi$-regularity of affine Schubert varieties that avoids the Mehta--Ramanathan criterion as used in \cite{Cas22,FHLR25}. Moreover, applying this criterion to wildly ramified groups in \cite{FHLR25} required additional casework, whereas our new proof is uniform across all groups, and it completes the last remaining open case of odd unitary groups when $p=2$. When $F$ has characteristic $0$, the hypothesis on the $\mathsf{P}$-action is satisfied for all $w$ in the $\mu$-admissible set of Kottwitz--Rapoport \cite{KR00} associated with some minuscule conjugacy class of geometric coweights $\mu$. It also holds when the simple reflections in $s_{\bullet}$ are distinct.

\begin{remark}
In \cite{Bha12}, Bhatt proved that Schubert varieties in the finite flag variety of $\mathrm{GL}_n$ in positive characteristic are derived splinters (an alternative name for globally $+$-regular), using inversion of adjunction. Bhatt asked in \cite[Remarks 7.8 and 7.10]{Bha12} if his methods could be generalized to general groups, and our proof of Theorem \ref{thm2} answers this question positively. Indeed, all classical Schubert varieties arise as particular affine Schubert varieties for $F$ of characteristic $p$, and their classical Demazure resolutions are based BSDH $k$-varieties.
\end{remark}

Let us now explain what would happen if we dropped the hypothesis on the $\mathsf{P}$-action on $\mathsf{S}_w$, which pertains solely to the case when $F$ has characteristic $0$. One can always construct a class of deperfections $\mathsf{S}_{s_{\bullet},q_\bullet}$, where $q_\bullet$ is a certain non-decreasing sequence of powers of $p$ defined inductively.
The factor $\mathsf{S}_{s_i,q_i}$ in the twisted product is the $\varphi_{q_i}$-twist of the canonical $\mathsf{B}(k)$-equivariant smooth deperfection $\mathsf{S}_{s_i,1}\simeq \bbP^1_k$ of $\mathsf{S}_{s_i}$. When the assumption on the $\mathsf{P}$-action is not satisfied, we are forced to twist every new factor to the right by a nonnegative power of $p$, which we have no control over and comes from trying to kill $p$-torsion of Witt rings. This ultimately prevents the theta divisor from being semi-ample, which hinders the inversion of adjunction.

\subsection{Local models}

Finally, we give an application to local models. Recall that \cite{AGLR22,GL24} proves the existence and uniqueness of normal flat $O_E$-schemes $\mathsf{M}_{\mu}$ with reduced special fiber representing a certain closed v-subsheaf of the Beilinson--Drinfeld Grassmannian $\mathrm{Gr}_{\mathcal{G}}$, provided either $\mu$ is minuscule or $F$ has characteristic $p$. In \cite{FHLR25} it was proved for all groups except wild odd unitary ones that the special fiber is moreover $\varphi$-split. Now, we can generalize this to all groups and prove it uniformly. This finishes the problem of determining the special fiber of $\mathsf{M}_{\mu}$ in full generality.

\begin{corollary}
	Assume $F$ has characteristic $p$ or $\mu$ is minuscule. Then, the special fiber of $\mathsf{M}_{\mu}$ equals the canonical deperfection $\mathsf{A}_{\mu}^{\mathrm{can}}$ in the sense of \cite{AGLR22} of the $\mu$-admissible locus. Moreover, $\mathsf{A}_{\mu}^{\mathrm{can}}$ is $\varphi$-split compatibly with every $\calG(O)$-stable closed subscheme.
\end{corollary}

Here, the canonical deperfection $\mathsf{A}_{\mu}^{\mathrm{can}}$ of the admissible locus is built out of gluing the $\mathsf{S}_w^{\mathrm{can}}$. We compute global sections of ample line bundles on the previous deperfection and the generic fiber of $\mathsf{M}_{\mu}$, cf.~the coherence conjecture of \cite{PR08}.
The $\varphi$-splitness yields higher vanishing of cohomology for ample line bundles, so we get an inclusion-exclusion type formula in terms of Schubert subvarieties: this is a combinatorial gadget thanks to the Demazure character formula, so we can check it on tame equicharacteristic $G$, already handled by Zhu \cite{Zhu14}. Similar ideas are used in \cite{Lou23} to prove normality of Schubert varieties embedded in $\mathsf{G}^{\mathrm{can}}/\mathsf{P}^{\mathrm{can}}$ when $F$ has equicharacteristic and $\pi_1(G)$ is $p$-torsion free. After this paper appeared, the work of He--Schremmer--Yu \cite{HSY26} settled Cohen--Macaulayness of $\mathsf{M}_\mu$ in full generality and uniformly so, building on the corollary above and their combinatorial methods.

\subsection{Acknowledgements}
We thank Johannes Anschütz, Sebastian Bartling, Bhargav Bhatt, Patrick Bieker, Ian Gleason, Eloísa Grifo, Arthur-César Le Bras, Xuhua He, Cédric Pépin, Thomas Polstra, Timo Richarz, Simon Riche, Peter Scholze, Karl Schwede, Austyn Simpson, Karen Smith, Felix Schremmer, Kevin Tucker, Thibaud van den Hove, Jakub Witaszek, and Jie Yang, and Qingchao Yu for helpful discussions. The authors used Claude and Gemini for literature review and grammar checking during the preparation of this article.
R.C.~was supported by the National Science
Foundation under Award No.~2103200 and No.~1840234. 
R.C.~also thanks Claremont McKenna College, the University of Michigan, and Sorbonne Paris Nord University for financial support during part of the preparation of this article.
J.L.~was supported by the Excellence Cluster of the Universität Münster, the ERC Consolidator Grant 770936 of Eva Viehmann, the DFG via the TRR 326 Geometry and Arithmetic of Uniformized Structures and the SFB 1442 Geometry: Deformations and Rigidity.

\section{Perverse \texorpdfstring{$\bbF_p$}{Fp}-sheaves and \texorpdfstring{$\varphi$}{phi}-singularities} Fix a prime number $p$.
For a scheme $X$ over $\bbF_p$, let $\varphi$ be the absolute Frobenius morphism. We will often be concerned with noetherian schemes which are $\varphi$-finite, meaning that $\varphi_* \mathcal{O}_X$ is a finite $\mathcal{O}_X$-module. By Kunz's theorem \cite[Theorem 2.5]{Kun76}, a $\varphi$-finite noetherian ring is excellent. Additionally, a noetherian $\varphi$-finite scheme admits a coherent dualizing complex \cite[Remark 13.6]{Gab04}. The proof in loc.~cit.~only applies when $X$ is affine, which is the only case we will use. Recall also that the perfection of a scheme $X$ is $X^{\mathrm{pf}}= \lim (\cdots \xrightarrow{\varphi} X \xrightarrow{\varphi} X)$. A deperfection of a perfect scheme $X$ is a scheme $X_0$ equipped with an isomorphism $X_0^{\mathrm{pf}} \cong X$. 

\subsection{Cartier modules}
Let $R$ be a ring over $\bbF_p$ and let $\varphi_* R$ be the $R$-module associated with $\varphi_*\calO_{\Spec (R)}$.
Recall from \cite{BB11} that a Cartier module over $R$ consists of an $R$-module $M$ with a map $\varphi_\ast M \to M$. Homomorphisms between Cartier modules must respect this map. A Cartier module $M$ is said to be nilpotent if $\varphi^e_\ast M \to M$ is zero for some $e \geq 0$. Furthermore, a Cartier module is said to be coherent if its underlying $R$-module is finite. We have the following decisive structure theorem for Cartier modules.

\begin{theorem}[Blickle--B\"ockle] \label{thm--CartierFinite}
    Let $R$ be a noetherian $\varphi$-finite ring and let $M$ be a coherent Cartier module. 
    \begin{enumerate}
        \item There exists a finite composition series $0 = M_0 \subset M_1 \subset \cdots \subset M_n = M$ by coherent Cartier submodules such that each $M_i/M_{i-1}$ is either nilpotent, or non-nilpotent and simple.
        \item If $M$ is a simple coherent Cartier module then $M$ has a unique associated prime $\frakp \in \Spec(R)$. Furthermore, $M \subset M_{\frakp}$, and the latter is a finite-dimensional vector space over $R/\frakp$.
    \end{enumerate}
\end{theorem}

\begin{proof}
    Part (1) is \cite[Proposition 4.23]{BB11}, and part (2) is proved in \cite[Propositions 4.14, 4.15]{BB11}.
\end{proof}

Important examples of coherent Cartier modules include the cohomology sheaves $\mathcal{H}^i(\omega_R^\bullet)$ of dualizing complexes on $\varphi$-finite noetherian  rings. Here the map $\varphi_\ast \mathcal{H}^i(\omega_R^\bullet) \to \mathcal{H}^i(\omega_R^\bullet)$ is obtained from exactness of $\varphi_*$ and the adjoint of the canonical 
isomorphism $\omega_R^\bullet \rightarrow \varphi^! \omega_R^\bullet$ from Grothendieck duality.

\subsection{$\varphi$-modules}
Let $R$ be an $\bbF_p$-algebra, and let $R[\varphi]$ be the non-commutative polynomial ring over $R$ in one variable, also denoted $\varphi$, subject to the relation $\varphi a = a^p \varphi$ for all $a \in R$. A left $R[\varphi]$-module is the same as an $R$-module $M$ with an $R$-linear map $M \rightarrow \varphi_*M$; note that the map goes in the direction opposite to that of Cartier modules. 

We recall a decisive structure result for $R[\varphi]$-modules closely related to Theorem \ref{thm--CartierFinite}. As in the case of Cartier modules, we say that an $R[\varphi]$-module $M$ is nilpotent if $M \rightarrow \varphi_*^e M$ is zero for some $e \geq 0$. Similarly, an $R[\varphi]$-module $M$ is co-finite if it is Artinian as an $R$-module. Important examples of co-finite $R[\varphi]$ modules include the local cohomology groups $H_{\mathfrak{m}}^i(R)$ of noetherian local $\bbF_p$-algebras $(R, \mathfrak{m})$, see \cite[Theorem 7.1.3]{BSh13}.

\begin{theorem}[Lyubeznik] \label{thm--varphi-modules}
    Let $(R, \mathfrak{m})$ be a noetherian local $\bbF_p$-algebra and let $M$ be a co-finite $R[\varphi]$-module. 
    \begin{enumerate}
        \item $M$ admits a finite composition series $0 = M_0 \subset M_1 \subset \cdots \subset M_n = M$ by co-finite $R[\varphi]$-submodules such that each $M_i/M_{i-1}$ is either nilpotent, or non-nilpotent and simple.
        \item The collection of non-nilpotent simple subquotients of $M$ is independent of the composition series.
    \end{enumerate}
\end{theorem}

\begin{proof}
See \cite[Theorem 4.7]{Lyu97}.
\end{proof}

It is worth mentioning the following special case of Lyubeznik's theorem, which has been proved via different means by various sets of authors.

\begin{corollary} \label{prop--cofinite}
    Let $(R, \mathfrak{m})$ be a noetherian local $\bbF_p$-algebra and let $M$ be a co-finite $R[\varphi]$-module.  Then some power of $\varphi$ annihilates $$\{a \in H_{\mathfrak{m}}^i(R) \: : \: \varphi^e(a) = 0 \text{ for some } e >0 \}.$$ 
\end{corollary}

\begin{proof}
    This follows from Theorem \ref{thm--varphi-modules}; see also \cite[Proposition 1.11]{HS77}, \cite[Proposition 4.4]{Lyu97}, \cite[Lemma 13.1]{Gab04} or \cite[Corollary 4.25]{BBL+23}.
\end{proof}

We conclude by explaining a precise relation between Cartier modules and $R[\varphi]$-modules. Suppose that $(R, \mathfrak{m})$ is a complete, local, noetherian and $\varphi$-finite $\bbF_p$-algebra. Following \cite[Tag 0A82]{StaProj}, we normalize the dualizing complex $\omega_R^\bullet$ so that $R\Gamma_{\mathfrak{m}}(\omega_R^\bullet) = E[0]$ lies in degree $0$, in which case $E$ is an injective hull of $R/\mathfrak{m}$.
Recall that Matlis duality $M \mapsto \Hom_R(M, E)$ gives an anti-equivalence between coherent and Artinian $R$-modules. Then Matlis duality also induces an anti-equivalence between coherent Cartier modules and co-finite $R[\varphi]$-modules \cite[Proposition 5.2]{BB11}. The integer $d:= \dim R$ is the largest integer such that $\mathcal{H}^{-d}(\omega_R^\bullet) \neq 0$ \cite[0AWN]{StaProj}; the cohomology sheaf $\omega_R := \mathcal{H}^{-d}(\omega_R^\bullet)$ is called the dualizing sheaf. For each $i$ there is a canonical isomorphism $\Hom_R(\mathcal{H}^{-i}(\omega_R^\bullet), E) \cong H_{\mathfrak{m}}^i(R)$  even if $R$ is not complete, e.g.~ see \cite[10.2.19]{BS13}, \cite[Tag 0AAK]{StaProj}.

\subsection{Perverse $\bbF_p$-sheaves}
For a scheme $X$ over $\bbF_p$ let $D(X, \bbF_p)$ be the derived category of \'etale $\bbF_p$-sheaves on $X$, and let $D_c^b(X, \bbF_p)$ be the bounded constructible subcategory. In this subsection we fix a perfect field $k$ of characteristic $p$. Every scheme of finite type over $k$ is automatically $\varphi$-finite.

\begin{definition}
Let $X$ be a $k$-scheme of finite type. For each point $x \in X$, fix a strict henselization $\calO_{x}^{\sh}$ of the local ring at $x$, and let $i_x \colon \overline{x} \rightarrow \Spec (\calO_{x}^{\sh})$ be the inclusion of the closed point.
We define the full subcategory $^{p}D^{\leq 0}(X, \mathbb{F}_p)$ (resp.~$^{p}D^{\geq 0}(X, \mathbb{F}_p)$) of $D(X, \bbF_p)$ consisting of $\mathcal{F}^\bullet \in D(X, \bbF_p)$ such that $\calH^n (i_x^* \mathcal{F}^\bullet) = 0$ for all $x \in X$ and $n > -\dim \overline{\{x\}}$ (resp.~$\mathcal{F}^\bullet$ has bounded below cohomology sheaves and $\calH^n (Ri_x^! \mathcal{F}^\bullet) = 0$ for all $x \in X$ and $n < -\dim \overline{\{x\}}$).  
\end{definition}

The following special case of a theorem of Gabber implies that the subcategories above give a t-structure on $D(X, \bbF_p)$, cf.~\cite[Theorem 11.5.4]{EK04a}. We call objects in the heart perverse $\mathbb{F}_p$-sheaves.

\begin{theorem}[Gabber] \label{thm-Gab04}
Let $X$ be a $k$-scheme of finite type.
\begin{enumerate}
\item The pair $({}^{p}D^{\leq 0}(X, \mathbb{F}_p), {}^{p}D^{\geq 0}(X, \mathbb{F}_p))$ gives rise to a t-structure on $D(X, \bbF_p)$.
\item The t-structure above restricts to a t-structure on $D_c^b(X, \bbF_p)$.
\item Every perverse subquotient of a constructible perverse $\bbF_p$-sheaf is constructible, i.e.~lies in $D_c^b(X, \bbF_p)$.
\item Every constructible perverse $\bbF_p$-sheaf has finite length.
\end{enumerate}
\end{theorem}

\begin{proof}
See \cite[Theorem 10.4, Corollary 12.4]{Gab04}.
\end{proof}

The t-structure on $D_c^b(X, \bbF_p)$ has also been studied in \cite{Cas22, BBL+23}.
By the topological invariance of the small \'etale site \cite[Tag 04DY]{StaProj}, we have a canonical equivalence $D_c^b(X, \bbF_p) \cong D_c^b(X^{\mathrm{pf}}, \bbF_p)$, so we also get a t-structure for perfections of $k$-schemes of finite type.

We now recall the notion of intermediate extension for perverse $\bbF_p$-sheaves. Let $\calF^\bullet$ be a constructible perverse $\bbF_p$-sheaf on $U$, and let $j \colon U \rightarrow X$ be an open immersion into a $k$-scheme of finite type. By taking perverse truncations of $Rj_*$ and $Rj_!$, we may define
$$j_{!*}\calF^\bullet:=  \text{Im} \: ({}^pj_! \calF^\bullet \rightarrow {}^pj_*\calF^\bullet).$$
Note that while ${}^pj_*\calF^\bullet$ may not be constructible, both ${}^pj_! \calF^\bullet$ and $j_{!*}\calF^\bullet$ are constructible. The intermediate extension $j_{!*}\calF^\bullet$ is characterized as the unique perverse extension of $\calF^\bullet$ with no quotients or subobjects supported on $X \setminus U$. If $i \colon X \setminus U \rightarrow X$ is the inclusion (with the reduced scheme structure), the latter conditions are equivalent to $i^*\calF^\bullet \in {}^{p}D^{\leq -1}(X \setminus U, \mathbb{F}_p)$ and $Ri^!(\calF^\bullet) \in {}^{p}D^{\geq 1}(X \setminus U, \mathbb{F}_p)$, respectively, by \cite[Lemma 2.7]{Cas22}.

\subsection{The Riemann--Hilbert correspondence}
We now recall the Riemann--Hilbert correspondence of Bhatt--Lurie \cite{BL19}. Let $R$ be an $\bbF_p$-algebra and let $(R, \varphi)$ be the ring $R$ regarded as an $R[\varphi]$-module via the Frobenius. For an $R$-algebra $S$, extension of scalars provides a functor from $R[\varphi]$-modules to $S[\varphi]$-modules, which is used implicitly in the following definition taken from \cite[Construction 2.3.1]{BL19}.
\begin{definition}
    Let $D(R[\varphi])$ be the derived category of $R[\varphi]$-modules. Define the functor $$\text{Sol}(-) := \underline{\text{RHom}}_{D(R[\varphi])} ((R, \varphi), -) \colon D(R[\varphi]) \to D(\Spec(R), \bbF_p).$$
\end{definition}

Informally, $\text{Sol}$ can be thought of as the derived functor of $\varphi$-invariants.
The functor $\text{Sol}$ is not an equivalence of categories because $D(R[\varphi])$ is too large. To solve this issue in the constructible case, Bhatt--Lurie define a notion of holonomicity, see \cite[Definition 4.1.1.]{BL19}. A holonomic $R[\varphi]$-module is an $R[\varphi]$-module isomorphic to one of the form $$M^\mathrm{pf} := \colim (M \rightarrow \varphi_*M \rightarrow \varphi_*^2 M \rightarrow \cdots )$$ for an $R[\varphi]$-module $M$ which is finite type as an $R$-module. Note that a holonomic $R[\varphi]$-module is in particular perfect, meaning that $M \rightarrow \varphi_*M$ is an isomorphism. Restriction of scalars along $R \rightarrow R^{\mathrm{pf}}$ identifies the categories of perfect $R[\varphi]$-modules and perfect $R^{\mathrm{pf}}[\varphi]$-modules by \cite[Proposition 3.4.3]{BL19}. By the following theorem, this is closely related to the topological invariance of the small \'etale site.

\begin{theorem}[Bhatt--Lurie] \label{thm--BhattLurie}
    Let $D_{\textnormal{hol}}(R[\varphi]) \subset D(R[\varphi])$ be the full subcategory of complexes with holonomic cohomology sheaves. Then $\text{Sol}$ restricts to an equivalence of categories $D_{\textnormal{hol}}(R[\varphi]) \cong D_c^b(\Spec(R), \bbF_p)$ which is t-exact for the standard t-structures on the source and target
\end{theorem}

\begin{proof}
See \cite[Theorem 7.4.1, Corollary 12.1.7]{BL19}.
\end{proof}

By t-exactness, if $M$ is a holonomic $R[\varphi]$-module then $\text{Sol}(M)$ is the \'etale sheaf on $\Spec(R)$ whose value on an \'etale $R$-algebra $S$ is
$$\text{Sol}(M)(S) = \{x \in M \otimes_R S \: : \: \varphi(x) = x\}.$$

\begin{remark} In \cite{BBL+23} the authors use a different definition of the perverse t-structure on $D_c^b(\Spec(R), \bbF_p)$, in terms of the Riemann--Hilbert correspondence and a perverse t-structure on coherent sheaves, but the two agree by \cite[Theorem 4.44]{BBL+23}. Correspondingly, our proofs of Theorem \ref{thm--CM-conditions} and Theorem \ref{prop-Fp-simple} below are quite different from their analogs \cite[Remark 4.40, Corollary 5.15]{BBL+23}. Our Theorem \ref{prop-Fp-simple} also differs from \cite[Corollary 5.15]{BBL+23} in that we allow a non-complete base and hence have to eliminate the possibility of branching behavior (with the help of \cite{DMP24}). 
\end{remark}

\subsection{Cohen--Macaulayness}
There are numerous equivalent definitions of Cohen--Macaulayness for a noetherian local ring, for example involving regular sequences, local cohomology, or a dualizing complex.
While there is no standard definition of Cohen--Macaulayness in the non-noetherian setting, the one in \cite[Definition 2.1]{Bha20} will be useful here (see also \cite[Remark 2.4]{Bha20}).

\begin{definition} \label{def-CM}
Let $X$ be a topologically noetherian scheme. We say that $X$ is Cohen--Macaulay if for every local ring $(R, \mathfrak{m})$ on $X$, the (Zariski) local cohomology groups $H_{\mathfrak{m}}^i(R):= R^i\Gamma_{\{\mathfrak{m}\}}(\mathcal{O}_{\Spec(R)})$ vanish for $i < \dim R$.
\end{definition}

If $R$ is a noetherian local $\bbF_p$-algebra then $H_{\mathfrak{m}}^i(R)$ has a canonical $R[\varphi]$-module structure as the cohomology of a Koszul complex of $R[\varphi]$-modules by \cite[Tag 0956]{StaProj}. In this case, $H_{\mathfrak{m}}^i(R)$ is finitely generated as an $R[\varphi]$-module, and since it is Artinian as an $R$-module it is even a module over the completion $\hat{R}$. Furthermore, if $M$ is an $R$-module ($R$ is still noetherian) we have $H^i_{\frakm}(M) = \colim_n \Ext^i_R(R/\frakm^n, M)$ by \cite[Tag 0955]{StaProj}.

\begin{lemma} \label{lemm--perfectloc}
    Let $(R, \mathfrak{m})$ be the perfection of a noetherian local $\bbF_p$-algebra $(R_0, \mathfrak{m}_0)$ of dimension $d$ with normalized dualizing complex $\omega_{R_0}^\bullet$. 
    \begin{enumerate}
        \item $H_{\mathfrak{m}}^i(R) = 0$ if and only if $H_{\mathfrak{m}_0}^i(R_0)$ is nilpotent.
        \item $H_{\mathfrak{m}}^{d}(R) \neq 0$.
        \item If $H_{\mathfrak{m}}^i(R) = 0$ for all $i < d$ then $R$ is equidimensional.
        \item If $R_0$ is $\varphi$-finite then $H^i_{\mathfrak{m}_0}(R_0)$ is nilpotent if and only if $\mathcal{H}^{-i}(\omega_{R_0}^\bullet)$ is nilpotent.
        \item If $R_0$ is $\varphi$-finite and $H_{\mathfrak{m}}^i(R) = 0$ for all $i < d$, then $\Spec(R)$ is Cohen--Macaulay in the sense of Definition \ref{def-CM}.
    \end{enumerate}
\end{lemma}

\begin{remark}
A noetherian local $\bbF_p$-algebra $(R_0, \frakm_0)$ such that $H_{\mathfrak{m}_0}^i(R_0)$ is nilpotent for $i < \dim R_0$ is called weakly $\varphi$-nilpotent in \cite{Mad19}, see also \cite{PQ19, Quy19}.
\end{remark}

\begin{proof}
    Since $R$ is an $R_0$-module and $\Spec(R) \cong \Spec(R_0)$ then $H_{\mathfrak{m}}^i(R) = H_{\mathfrak{m}_0}^i(R)$. Next, $H_{\mathfrak{m}_0}^i(R) = \colim_n \Ext^i_{R_0}(R_0/\frakm_0^n, \colim_e \varphi^e_* R_0)$. By taking a resolution of $R_0/\frakm_0^n$ by finite free $R_0$-modules and using that filtered colimits are exact \cite[Tag 00DB]{StaProj}, the inner colimit over $e$ commutes with $\Ext^i_{R_0}(R_0/\frakm_0^n, -)$. Then by exchanging the colimits and using exactness of $\varphi_*$ to commute the latter with  $H_{\mathfrak{m}_0}^i(-)$, we get
    \begin{equation} \label{eqn-local-coh-Frob} H_{\mathfrak{m}}^i(R) = \underset{e}{\colim} \, \varphi_*^e H_{\mathfrak{m}_0}^i(R_0) = H_{\mathfrak{m}_0}^i(R_0)^{\mathrm{pf}}.
\end{equation}
Thus, $H_{\mathfrak{m}}^i(R) = 0$ if and only if every element of $H_{\mathfrak{m}_0}^i(R_0)$ is annihilated by some power of $\varphi$. Now (1) follows from Corollary \ref{prop--cofinite}. 

For (2), we consider two cases. If $d = 0$ then $H_{\mathfrak{m}}^{0}(R) = R$ is nonzero. On the other hand, if $d > 0$, then by (1) we need only show that $H_{\mathfrak{m}_0}^d(R_0)$ is not nilpotent. But if $H_{\mathfrak{m}_0}^d(R_0)$ were nilpotent, then since it is finitely generated as an $R_0[\varphi]$-module, it would also be finitely generated as an $R_0$-module, which is impossible by \cite[Corollary 7.3.3]{BSh13}.

Part (3) follows from the proof of \cite[Proposition 2.8(3)]{PQ19}; we reproduce the argument here for completeness.
If $R$ is not equidimensional, let $\mathfrak{p} \subset R_0$ be a minimal prime such that $n: = \dim R_0/\mathfrak{p} < d = \dim R_0$, and let $I$ be the intersection of the other minimal primes.
Then we have an exact sequence of $R_0[\varphi]$-modules
$$0 \rightarrow R_0 \rightarrow R_0/\mathfrak{p} \oplus R_0/I \rightarrow R_0/(\mathfrak{p} + I) \rightarrow 0$$ where $\dim R_0/(\mathfrak{p} + I) < n$. Applying $R\Gamma_{\{\mathfrak{m}_0\}}$ gives a surjection $H_{\frakm_0}^n(R_0) \rightarrow H^n_{\frakm_0}(R_0/\mathfrak{p})$ by \cite[Tag 0DXC]{StaProj}. This implies that $H^n_{\frakm_0}(R_0/\mathfrak{p})$ is nilpotent, which contradicts part (2). 

Next, we note that if $R_0$ is complete then (4) follows immediately from Matlis duality. Remarkably this statement is true even if $R_0$ is not complete, as was observed in \cite[Lemma 2.3]{ST17}. The argument is similar to \cite[Proposition 4.3]{Sch09}, using that the double Matlis duality functor is isomorphic to $(-) \otimes_{R_0} \hat{R}_0$ on finite $R_0$-modules, and faithful flatness of $R_0 \rightarrow \hat{R}_0$; we refer to loc.~cit for more details. 

For each prime $\frakp \subset R_0$, the localization $(\omega_{R_0}^\bullet)_{\frakp}$ is a dualizing complex for $(R_0)_{\frakp}$, and since $R_0$ is equidimensional, $(\omega_{R_0})_{\frakp}$ is a dualizing sheaf for $(R_0)_{\frakp}$ by \cite[Proposition 2.3.2]{Smi93}. Thus, (5) follows from (4).
\end{proof}

Since  $H_{\mathfrak{m}}^i(R) = H_{\mathfrak{m}_0}^i(R)$, then $H_{\mathfrak{m}}^i(R)$ has a canonical $R[\varphi]$-module structure as the cohomology of a Koszul complex, constructed from finitely many generators of $\frakm$ up to radical, independent of the chosen deperfection. 
The canonicity of the $R[\varphi]$-module structure can also be deduced by noting that the action of $\varphi$ comes from applying $H_{\mathfrak{m}}^i(-)$ to the Frobenius map $R \to \varphi_* R$. 
Furthermore, $H_{\mathfrak{m}}^i(R)$ is a perfect $R[\varphi]$-module in the sense of \cite[Definition 3.2.1]{BL19}.

We now introduce some notation for working with perfections. If $R$ is an $\bbF_p$-algebra and $r \in R$, we denote by $r^{1/p^{e}} \in R^{\mathrm{pf}}$ the $p^{e}$-th root. If $M$ is an $R[\varphi]$-module and $m \in M$, we denote by $\varphi^{-e}(m) \in M^{\mathrm{pf}}$ the image of $m$ under the map $\varphi_*^eM \to M^{\mathrm{pf}}$. We use this notation in the proof of the following structure result for perfect $\varphi$-modules.

\begin{proposition}
    \label{thm--perfectstructure}
    Let $(R, \mathfrak{m})$ be the perfection of a noetherian local $\bbF_p$-algebra $(R_0, \mathfrak{m}_0)$.  
    \begin{enumerate}
        \item The functor $M_0 \mapsto M_0^\mathrm{pf}$ from $R_0[\varphi]$-modules to $R[\varphi]$-modules is exact.
        \item If $M_0$ is a co-finite $R_0[\varphi]$-module then $M_0^\mathrm{pf}$ has finite length as an $R[\varphi]$-module.
        \item If $M_0$ is a co-finite, non-nilpotent and simple $R_0[\varphi]$-module then $M_0^\mathrm{pf}$ is a simple, nonzero $R[\varphi]$-module.
    \end{enumerate}
\end{proposition}

\begin{proof}
    Part (1) follows from the exactness of $\varphi_*$. For (2), we take the perfection of a composition series for $M_0$ as in Theorem \ref{thm--varphi-modules}. This kills the nilpotent subquotients, so (2) will then follow from (3). For (3), let $M:= M_0^\mathrm{pf}$, and let $m \in M$ be nonzero. We need to show that $R[\varphi] \cdot m = M$. Write $m = \varphi^{-e}(m')$ for some $m' \in M_0$ and $e \geq 0$. It suffices to show that for all $n \in M_0$ and $f \geq 0$, we have $\varphi^{-e-f}(n) \in R[\varphi] \cdot m$. By simplicity $R_0[\varphi] \cdot \varphi^f(m') = M_0$, so there exist $r_i \in R_0$ such that $\sum_i r_i \cdot \varphi^{i+f} (m') = n$. Then we conclude since $\sum_i r_i^{1/p^{e+f}} \cdot \varphi^i (m) = \varphi^{-e-f}(n).$
\end{proof}

The following result is the key input in the proof of Theorem \ref{thm--CM-conditions} below.

\begin{proposition} \label{lemm--punctureCM}
    Let $(R, \frakm)$ be the perfection of a $\varphi$-finite noetherian local $\bbF_p$-algebra $(R_0, \frakm_0)$. Suppose that $R$ is equidimensional, the punctured spectrum of $(R, \frakm)$ is Cohen--Macaulay, and $R/\frakm$ is algebraically closed. Then $\Spec(R)$ is Cohen--Macaulay if and only if for all $i < \dim R$, there does not exist a nonzero element $x \in H^i_{\frakm}(R)$ such that $\varphi(x) =x$.
\end{proposition}

\begin{proof}
    The necessity of the condition on $\varphi$-fixed elements is clear. For sufficiency,  by Lemma \ref{lemm--perfectloc} we may suppose for contradiction that $\mathcal{H}^{-i}(\omega_{R_0}^\bullet)$ has a non-nilpotent simple Cartier subquotient $M$ for some $i \neq \dim R$. By Theorem \ref{thm--CartierFinite} and our assumption on the punctured spectrum, the unique associated prime of $M$ must be $\frakm_0$, so $M$ has finite length as an $R_0$-module. Since $M$ was arbitrary, $\mathcal{H}^{-i}(\omega_{R_0}^\bullet)$ has finite length as an $R_0$-module up to nilpotents. Let $E$ be an injective hull of $R_0/\frakm_0$. Then $\Hom_{R_0}(\mathcal{H}^{-i}(\omega_{R_0}^\bullet), E) \cong H^i_{\frakm_0}(R_0)$, even as $R_0[\varphi]$-modules by \cite[Lemma 5.1]{BB11}, so $H^i_{\frakm_0}(R_0)$ has finite length as an $R_0$-module up to nilpotents.

    Now by Proposition \ref{thm--perfectstructure}, $H^i_{\frakm}(R)$ has a finite composition series with holonomic subquotients in the sense of \cite[Definition 4.1.1]{BL19}. Since the abelian category of holonomic $R_0[\varphi]$-modules is closed under extensions, see \cite[Corollary 4.3.3, Remark 3.2.2]{BL19}, then $H^i_{\frakm}(R)$ is a holonomic $R_0$-module which is moreover set-theoretically supported on $\Spec(R/\frakm)$ (holonomicity can also be deduced from \cite[Lemma 2.16]{Bha20}). By construction the functor $\text{Sol}$ is compatible with pullback, so that  $\text{Sol}(H^i_{\frakm}(R))$ is an \'etale sheaf supported on $\Spec(R/\frakm)$. Since $R/\frakm$ is algebraically closed, an \'etale sheaf vanishes if and only if its global sections vanish. Now we conclude using that $\text{Sol}$ is an equivalence when restricted to holonomic modules by Theorem \ref{thm--BhattLurie}.
\end{proof}

\begin{lemma} \label{lemm--catenary}
    Let $(R, \frakm)$ be a noetherian local ring which is equidimensional and universally catenary. Let $R \rightarrow S$ be an \'etale ring map and let $\frakq$ be a prime of $S$ lying above $\frakm$. Then the localization $S_\frakq$ is equidimensional and $\dim S_{\frakq} = \dim R$.
\end{lemma}

\begin{proof}
Let $\frakq_0$ be a minimal prime of $S$ contained in $\frakq$. Then $\frakp_0 := R \cap \frakq_0$ is minimal by flatness, and $R/\frakp_0$ is universally catenary by \cite[Tag 00NK]{StaProj}. By the dimension formula in \cite[Tag 02IJ]{StaProj} applied to the ring extension $R/\mathfrak{p}_0 \subset S/\frakq_0$, we have $\text{ht}(\frakm / \frakp_0) = \text{ht}(\frakq/\frakq_0)$.  
\end{proof}

\begin{lemma} \label{lemm--invariance}
    Let $(R, \frakm)$ be a local $\mathbb{F}_p$-algebra. Then if $R$ satisfies any of the following three properties, so does the strict henselization $(R^{\sh}, \frakm^{\sh})$.
    \begin{enumerate}
        \item $R$ is noetherian.
        \item $R$ is $\varphi$-finite.
        \item $R$ is excellent and equidimensional.
    \end{enumerate}
\end{lemma}

\begin{proof}
    Property (1) is part of \cite[Proposition 18.8.8]{EGAIV4}. Property (2) is in the proof of \cite[Theorem 4.1]{BCGST19}, the point being that $R^{\sh}$ is an ind-\'etale $R$ algebra, so $\varphi_* R \otimes_{R} R^{\sh} = \varphi_* R^{\sh}$. For property (3), excellence is preserved by the last remark in \cite[Ch. 1 \S 1]{FK88}, and it remains to prove equidimensionality.

    Let $\frakq_0$ be a minimal prime of $R^{\sh}$. Then $\frakp_0 := R \cap \frakq_0$ is minimal by flatness. As in \cite[Tag 06LK]{StaProj}, write $R^{\sh} = \colim_i R_i$ as a direct limit of local rings $R_i$ which are localizations of \'etale $R$-algebras faithfully flat over $R$. Then $R_i \rightarrow R^{\sh}$ is faithfully flat for all $i$ by \cite[Tag 00U7, Tag 05UT]{StaProj}. We have the minimal primes $\frakp_i := R_i \cap \frakq_0$ in the $R_i$. 
    By prime avoidance, for all large enough $i$, $\frakq_0$ is the only minimal prime of $R^{\sh}$ lying over $\frakp_i$. Fix such an $i$ and let $C=(\frakQ_0 \subset \cdots \subset \frakQ_n)$ be a maximal chain of primes in $R^{\sh}/\frakp_i R^{\sh}$. By faithful flatness of $R_i/{\frakp_i} \rightarrow R^{\sh}/\frakp_i$ and the going down property, we have $n \geq \dim R_i/\frakp_i$. But by Lemma \ref{lemm--catenary}, $\dim R_i/\frakp_i = \dim R$. Since $\dim R = \dim R^{\sh}$ by \cite[Tag 06LK]{StaProj}, this implies that the preimage $C \cap R^{\sh}$ is also a maximal chain of primes in $R^{\sh}$. By our choice of $i$, the minimal prime of $C \cap R^{\sh}$ is $\frakq_0$. Thus, for an arbitrary minimal prime $\frakq_0$ we have exhibited a chain $C$ of primes starting at $\frakq_0$ and of length $\dim R^{\sh}$, so $R^{\sh}$ is equidimensional.
\end{proof}

\begin{lemma} \label{lemm--loccohstrict}
     Let $(R, \mathfrak{m})$ be the perfection of a noetherian local $\bbF_p$-algebra, and let $(R^{\sh}, \mathfrak{m}^{\sh})$ be the strict henselization.
     \begin{enumerate}
         \item $H^i_{\mathfrak{m}^{\sh}}(R^{\sh}) = H^i_\mathfrak{m}(R) \otimes_R R^{\sh}.$
         \item $H^i_\mathfrak{m}(R) = 0$ if and only if $H^i_{\mathfrak{m}^{\sh}}(R^{\sh}) = 0$. 
         \item If the punctured spectrum of $(R, \mathfrak{m})$ is Cohen--Macaulay, so is the punctured spectrum of  $(R^{\sh}, \mathfrak{m}^{\sh})$.  
     \end{enumerate}
\end{lemma}

\begin{proof}
    It follows from topological invariance of the small \'etale site \cite[Tag 04DY]{StaProj} that perfection commutes with strict henselization, so that $R^{\sh}$ is topologically noetherian.
    Since $\mathfrak{m}^{\sh} = \frakm R^{\sh}$ then $H^i_{\frakm^{\sh}}(R^{\sh}) = H^i_{\frakm}(R^{\sh})$. Via the description of $H^i_{\frakm}(R)$ as the cohomology of a Koszul complex of $R$-modules, (1) and (2) follow from faithful flatness of $R \rightarrow R^{\sh}$, see \cite[Tag 07QM]{StaProj}. For (3), let $\mathfrak{q}$ be a non-maximal prime of $R^{\sh}$, and let $\mathfrak{p}$ be its preimage in $R$. We claim that $\frakq R^{\sh}_{\frakq}= \frakp R^{\sh}_{\frakq}$ and $\dim R^{\sh}_{\frakq} = \dim R_{\frakp}$. Granting these claims, $H^i_\frakq(R_\frakq^{\sh}) = H^i_\frakp(R_\frakp) \otimes R^{\sh}_\mathfrak{q}$ and (3) follows. 
    
    To prove the claims, use \cite[Tag 06LK]{StaProj} to write $R^{\sh} = \colim_i R_i$ as a direct limit of local rings $R_i$ which are localizations of \'etale $R$-algebras. If $\mathfrak{p}_i$ is the preimage of $\mathfrak{q}$ in $R_i$, then $R^{\sh}_{\frakq} = \colim_i (R_i)_{\frakp_i}$. We have $\frakp_i (R_i)_{\frakp_i} = \frakp (R_i)_{\frakp_i}$ by \cite[Tag 00U4]{StaProj}, so the first claim follows. For the claim about dimensions, note that $\dim R^{\sh}_{\frakq} \geq \dim R_{\frakp}$ by faithful flatness and the going down property. For the other direction, if $\mathfrak{q}_0 \subset \mathfrak{q}_1 \subset \cdots \subset \mathfrak{q}_n= \mathfrak{q}$ is a chain of primes in $R^{\sh}$, then for some $i$ this restricts to a chain of primes in $R_i$ of the same length. Now we conclude since $\dim (R_i)_{\frakp_i} = \dim R_{\frakp}$ by \cite[Tag 07QP]{StaProj}.
\end{proof}

We now come to the main result of this subsection. In the proof, we use the Artin--Schreier sequence
\begin{equation} \label{eq--AS}
 0\rightarrow \bbF_p \rightarrow \mathcal{O}_{\Spec(R)} \xrightarrow{\varphi - 1} \mathcal{O}_{\Spec(R)} \rightarrow 0
 \end{equation}
to translate between quasi-coherent cohomology and \'etale cohomology. This sequence is exact in the \'etale topology on $\Spec(R)$ for any $\bbF_p$-algebra $R$. If $(R, \frakm)$ is the perfection of a local $\bbF_p$-algebra $R_0$ then $\mathcal{O}_{\Spec(R)}$
may be viewed as a quasi-coherent sheaf on $\Spec(R_0)$, and then \eqref{eq--AS} is also exact in the \'etale topology on $\Spec(R_0)$. In particular, if $R_0$ is noetherian then $R^i\Gamma_{\{\mathfrak{m}\}}(\mathcal{O}_{\Spec(R)}) = H^i_{\frakm_0}(R)$ is the same in both the Zariski and \'etale topology by \cite[Tag 04DY]{StaProj} and descent for quasi-coherent sheaves.

\begin{theorem} \label{thm--CM-conditions}
Let $k$ be a perfect field of characteristic $p$ and let $X$ be a scheme isomorphic to the perfection of a connected finite type $k$-scheme. Then the following are equivalent.
\begin{enumerate}
    \item $X$ is Cohen--Macaulay in the sense of Definition \ref{def-CM}.
    \item The shifted constant sheaf $\bbF_p[\dim X] \in D_c^b(X, \mathbb{F}_p)$ is perverse.
\end{enumerate}
Furthermore, $X$ is equidimensional in both cases.
\end{theorem}

\begin{proof}
We first show that (1) implies (2). By Lemma \ref{lemm--perfectloc}, $X$ is equidimensional. It is immediate that $\bbF_p[\dim X] \in {}^{p}D^{\leq 0}(X, \mathbb{F}_p)$. To prove that $\bbF_p[\dim X] \in {}^{p}D^{\geq 0}(X, \mathbb{F}_p)$, fix a point $x \in X$. Let $(R, \mathfrak{m})$ be the strict henselization of the corresponding perfect local ring. Let $d:= \dim R$, which also agrees with the dimension before strict henselization by \cite[Tag 06LK]{StaProj}. By Lemma \ref{lemm--loccohstrict}, (1) is equivalent to the statement that $H^i_{\frakm}(R) = 0$ for all $i < d$ and points $x \in X$.
On the other hand, as $X$ is equidimensional, (2) is equivalent to the statement that for all points $x$, we have $\calH^i(Ri_x^! \bbF_p) = 0$ for $i < d$, where $i_x \colon \overline{x} \rightarrow \Spec(R)$ is the inclusion of the closed point. Since $\overline{x}$ is a geometric point, $\calH^i(Ri_x^! \bbF_p) = R^i\Gamma_{\{ \mathfrak{m}\}}(\bbF_p)$. Now the fact that (1) implies (2) follows from applying $R\Gamma_{\{ \mathfrak{m}\}}$ to the Artin--Schreier sequence \eqref{eq--AS}.

Next, we show that (2) implies $X$ is equidimensional. A straightforward argument using the Artin--Schreier sequence  as above shows that if $Y$ is Cohen--Macaulay, irreducible, and of finite type over $k$, then $\bbF_p[\dim Y] \in D_c^b(Y, \bbF_p)$ is perverse. Now let $X_0$ be a deperfection of $X$ by a finite type $k$-scheme, which we may assume is reduced. 
Since $X_0$ is of finite type over the perfect field $k$, we may let $Y$ be a smooth dense open subscheme of $X_0$, so the perversity of $\bbF_p[\dim Y]$ implies that $X_0$ is equidimensional.

To show that (2) implies (1), we proceed by descending induction on $\dim \overline{\{x\}}$. By Lemma \ref{lemm--loccohstrict}, it suffices to show the strict henselization  $(R, \mathfrak{m})$ of the local ring at $x$ satisfies $H^i_{\mathfrak{m}}(R) = 0$ for $i < d:= \dim R = \dim X - \dim \overline{\{x\}}$, and we may assume the punctured spectrum of $(R, \mathfrak{m})$ is Cohen--Macaulay.
By Lemma \ref{lemm--invariance} the hypotheses of Proposition \ref{lemm--punctureCM} are satisfied, so we are reduced to checking that $H^i_{\mathfrak{m}}(R)$ has trivial $\varphi$-invariants for $i < d$.
But this condition follows from induction and the long exact sequence obtained from applying $R\Gamma_{\{ \mathfrak{m}\}}$ to the Artin--Schreier sequence \eqref{eq--AS}, together with the perversity of $\bbF_p[\dim X]$.
\end{proof}

\subsection{$\varphi$-rationality}
Next we discuss a perfect notion of $\varphi$-rationality. Recall that by a theorem of Smith \cite{Smi97}, an excellent local $\bbF_p$-algebra $(R_0, \mathfrak{m}_0)$ is $\varphi$-rational if and only if it is Cohen--Macaulay and $H^{\dim R_0}_{\frakm_0}(R_0)$ is a simple $R_0[\varphi]$-module (this differs from the original definition in terms of tight closure, cf.~\cite{FW89, HH94a}). The first author showed in \cite[Theorem 1.7]{Cas22} that if $X$ is an irreducible scheme of finite type over an algebraically closed field, all of whose local rings are $\varphi$-rational, then $\bbF_p[\dim X]$ is simple as a perverse sheaf. 
In the opposite direction, we will show that if $\bbF_p[\dim X]$ is simple, then $X^{\mathrm{pf}}$ is Cohen--Macaulay and $H^{\dim R_0}_{\frakm_0}(R_0)^{\mathrm{pf}}$ is simple. The latter properties are encapsulated by the following definition.

\begin{definition}
    Let $X$ be the perfection of a noetherian $\varphi$-finite $\bbF_p$-scheme. We say that $X$ is $\varphi$-rational if it is Cohen--Macaulay in the sense of Definition \ref{def-CM}, and, for every local ring $(R, \frakm)$ on $X$, the top local cohomology group $H_{\frakm}^{\dim R}(R)$ is a simple $R[\varphi]$-module.
\end{definition}

It will be useful to characterize $\varphi$-rationality of a perfect scheme in terms of a property of one (equivalently, every) deperfection. The latter property turns out to be $\varphi$-nilpotence, first introduced by Blickle--Bondu \cite{BB05} under the name close to $F$-rational, and further studied e.g.~in \cite{ST17, PQ19, DMP24, KMPS23}.

\begin{definition} \label{def--rational}
Let $(R, \frakm)$ be a $\varphi$-finite noetherian local $\bbF_p$-algebra of dimension $d$. 
\begin{enumerate}
\item The tight closure of the zero submodule in $H_{\frakm}^d(R)$, denoted $0_{H_{\frakm}^d(R)}^*$, is the $R[\varphi]$-submodule consisting of all elements $x \in H_{\frakm}^d(R)$ such that there exists some $c \in R$ not contained in any minimal prime with the property that $c\varphi^e(x) = 0$ for all $e \gg 0$. 
\item The ring $R$ is said to be $\varphi$-nilpotent if each of the $R[\varphi]$-modules $$H^0_\frakm(R), \ldots, H^{d-1}_{\frakm}(R), 0_{H_{\frakm}^d(R)}^*$$ is nilpotent. 
\end{enumerate}
\end{definition}

By \cite[Proposition 2.8 (2)]{PQ19}, $R$ is $\varphi$-nilpotent if and only if its reduction is $\varphi$-nilpotent, so we will usually assume $R$ is reduced. 
To relate this  notion to coherent objects, note that
Matlis duality gives a canonical pairing $f \colon H_{\frakm}^d(R) \otimes_{\hat{R}} \omega_{\hat{R}} \rightarrow E$. The parameter test module $\tau(\omega_R) \subset \omega_R$ is the Cartier submodule consisting of all $\eta \in \omega_R$ such that $f(x \otimes \eta) = 0$ for all $x \in 0_{H_{\frakm}^d(R)}^*$. The parameter test module is well-behaved under localization by \cite[Proposition 3.1]{HT04} (cf.~\cite[Proposition 3.2 (e)]{Bli13}), completion by \cite[Proposition 3.2]{HT04}, and more generally under flat base change when the residue field extension is separable, thanks to \cite[Lemma 1.5]{ST17}. By construction, $0_{H_{\frakm}^d(R)}^*$ is the Matlis dual of $\omega_R/ \tau(\omega_R)$, even if $R$ is not complete. When combined with Lemma \ref{lemm--perfectloc}, the following gives a characterization of $\varphi$-nilpotence in terms of $\omega_R^\bullet$.

\begin{lemma} \label{lemm--dual}
Let $(R, \frakm)$ be a reduced, $\varphi$-finite noetherian local $\bbF_p$-algebra of dimension $d$. Then $0_{H_{\frakm}^d(R)}^*$ is nilpotent if and only if $\omega_R/ \tau(\omega_R)$ is nilpotent.
\end{lemma}

\begin{proof} 
It is observed \cite[Lemma 2.3]{ST17} that this follows from an argument similar to \cite[Lemma 2.1]{HT04}.
\end{proof}

\begin{lemma} \label{lemm--rationalequiv}
Let $(R, \frakm)$ be the perfection of a $\varphi$-finite noetherian local $\bbF_p$-algebra $(R_0, \frakm_0)$ of dimension $d$. Then the following are equivalent.
\begin{enumerate}
\item $\Spec(R)$ is $\varphi$-rational in the sense of Definition \ref{def--rational}.
\item $\Spec(R^{\sh})$ is $\varphi$-rational in the sense of Definition \ref{def--rational}.
\item $R_0$ is $\varphi$-nilpotent.
\item $R_0^{\sh}$ is $\varphi$-nilpotent.
\end{enumerate}
Furthermore, if any of the above conditions is satisfied then $R$ is geometrically unibranch. 
\end{lemma}

\begin{proof} By \cite[Proposition 2.8]{PQ19} we may assume that $R_0$ is reduced.
The equivalence of (3) and (4) then follows from Lemma \ref{lemm--loccohstrict} and Lemma \ref{lemm--dual}, together with the compatibility of $\tau(\omega_{R_0})$ and its structure map as a Cartier module under faithfully flat base change to $R_0^{\sh}$ as in \cite[Proposition 2.4 (4)]{ST17}, cf.~\cite[Theorem 4.4]{KMPS23}. Once we prove the equivalence of (1) and (3), the equivalence with (2) will then follow.

First suppose that $R_0$ is $\varphi$-nilpotent. Then the completion $\hat{R}_0$ is also $\varphi$-nilpotent by \cite[Proposition 2.8 (4)]{PQ19}. Thus $\hat{R}_0$ is a domain by \cite[Theorem 3.1]{DMP24}, and after the equivalence of (1)-(4) is established, loc.~cit.~will also imply the final claim about geometric unibranchedness. By \cite[Theorem 3.1.4]{Smi93}, $0_{H_{\frakm_0}^d(R_0)}^*$ is the unique maximal proper $R_0[\varphi]$-submodule of $H_{\frakm_0}^d(R_0)$. Since this submodule is nilpotent, then Proposition \ref{thm--perfectstructure} implies that $H_{\frakm}^{d}(R)$ is a simple $R[\varphi]$-module. Furthermore, $H_{\frakm}^{i}(R) = 0$ for $i < d$ by Lemma \ref{lemm--perfectloc}. For every prime $\frakp \subset R_0$ the localization $(R_0)_{\frakp}$ is $\varphi$-nilpotent by \cite[Proposition 2.4 (3)]{ST17} or ~\cite[Corollary 5.17]{PQ19}, so the same arguments apply to $(R_0)_{\frakp}$ and hence $\Spec(R)$ is $\varphi$-rational.

Now suppose that $\Spec(R)$ is $\varphi$-rational. Then $H^i_{\frakm_0}(R_0)$ is nilpotent for $i < d$ and $R_0$ is equidimensional by Lemma \ref{lemm--perfectloc}. By \cite[Corollary 3.9]{Bli04} (and the surrounding discussion if $R_0$ is not complete), $0_{H_{\frakm_0}^d(R_0)}^*$ is the intersection of the maximal proper $R_0[\varphi]$-submodules of $H_{\frakm_0}^d(R_0)$. Furthermore, the quotient of $H_{\frakm_0}^d(R_0)$ by each of these maximal proper $R_0[\varphi]$-submodules is non-nilpotent by \cite[Theorem 3.8]{Bli04}. Since $H_{\frakm_0}^d(R_0)$ is simple up to nilpotents then  $0_{H_{\frakm_0}^d(R_0)}^*$ must be nilpotent, so $R_0$ is $\varphi$-nilpotent.
\end{proof}

\begin{lemma} \label{lemm--rationalpuncture}
    Let $(R, \frakm)$ be the perfection of a $\varphi$-finite noetherian local $\bbF_p$-algebra $(R_0, \frakm_0)$. If the punctured spectrum of $(R, \frakm)$ is $\varphi$-rational, then so is the punctured spectrum of  $(R^{\sh}, \frakm^{\sh})$. 
\end{lemma}

\begin{proof}
   By Lemma \ref{lemm--rationalequiv} the punctured spectrum of $(R_0, \frakm_0)$ is $\varphi$-nilpotent, and it suffices to show the same is true of $(R_0^{\sh}, \frakm_0^{\sh})$. We can assume $R_0$ is reduced. If $\frakq \subset R_0^{\sh}$ is a non-maximal prime then it lies over some non-maximal prime $\frakp \in R_0$.  The map of local rings $((R_0)_{\frakp}, \frakp (R_0)_{\frakp}) \rightarrow ((R_0^{\sh})_\frakq, \frakq(R_0^{\sh})_\frakq)$ is faithfully flat and the residue field extension is separable (for separability, use that $R^{\sh}_0$ is a filtered colimit of \'etale $R_0$-algebras as in \cite[Tag 04GW]{StaProj} and apply \cite[Tag 00U4]{StaProj}). Since $(R_0)_{\frakp}$ is $\varphi$-nilpotent, so is $(R_0^{\sh})_\frakq$ by \cite[Proposition 2.4 (4)]{ST17}.
\end{proof}

The following result generalizes \cite[Proposition 2.5]{ST17} to the case where the deperfected punctured spectrum of $(R_0, \frakm_0)$ is $\varphi$-nilpotent instead of $\varphi$-rational.

\begin{proposition} \label{lemm--puncturerational}   Let $(R, \frakm)$ be the perfection of a $\varphi$-finite noetherian local $\bbF_p$-algebra $(R_0, \frakm_0)$ of dimension $d > 0$. Suppose that $R$ is equidimensional, the punctured spectrum of $(R, \frakm)$ is $\varphi$-rational, and $R/\frakm$ is algebraically closed. Then $\Spec(R)$ is $\varphi$-rational if and only if for all $i$, there does not exist a nonzero element $x \in H^i_{\frakm}(R)$ such that $\varphi(x) =x$.
\end{proposition}

\begin{proof}
First suppose that $\Spec(R)$ is $\varphi$-rational. We may assume that $R_0$ is reduced. By Lemma \ref{lemm--punctureCM} we only need to deal with the conditions on $H^{d}_{\frakm}(R)$. As in the proof of Lemma \ref{lemm--rationalequiv}, $R_0$ is a domain and $H^{d}_{\frakm_0}(R_0)$ has a unique simple $R_0[\varphi]$-module quotient $M_0$, which is also non-nilpotent. We claim that  $\text{Ann}_{R_0}(M_0) = (0)$. This follows from $\text{Ann}_{\hat{R}_0}(M_0) = (0)$, which in turn follows from the fact that Matlis duality preserves annihilators by \cite[10.2.14]{BS13} and torsion-freeness of the dualizing sheaf $\omega_{\hat{R}_0}$, see \cite[Tag 0AWK]{StaProj}. Since $R$ is $\varphi$-rational then $H_{\frakm}^d(R)= M_0^{\mathrm{pf}}$, and it follows that $\text{Ann}_{R}(H_{\frakm}^d(R)) = (0)$.
Now if there exists a nonzero $x \in H_{\frakm}^d(R)$ with $\varphi(x) = x$, then by simplicity $x$ generates $H_{\frakm}^d(R)$ as an $R$-module. But every element of $H_{\frakm}^d(R)$ is also annihilated by some collection of elements which generate $\frakm$ up to radical, so $\frakm = (0)$, a contradiction since $d > 0$.

Now suppose the punctured spectrum of $(R_0, \frakm_0)$ is $\varphi$-nilpotent. Again, we only need to deal with the conditions on $H^{d}_{\frakm}(R)$, and $R_0$ is equidimensional by Lemma \ref{lemm--perfectloc}. For contradiction we may assume that $R_0$ is reduced and $\omega_{R_0}/\tau(\omega_{R_0})$ is non-nilpotent (Lemma \ref{lemm--dual}). Let $M$ be a simple non-nilpotent Cartier subquotient of $\omega_{R_0}/\tau(\omega_{R_0})$, and let $\frakp \subset R_0$ be its unique associated prime (Theorem \ref{thm--CartierFinite}). By the compatibility of $\tau(\omega_{R_0})$ with localization, see \cite[Proposition 3.1]{HT04}, our assumption on the punctured spectrum implies $\frakp = \frakm_0$. Now we conclude as in the proof of Proposition \ref{lemm--punctureCM}. Briefly, the Matlis dual of $\omega_{R_0}/\tau(\omega_{R_0})$ is $0^*_{H_{\frakm_0}^d(R_0)}$, which therefore has finite length as an $R_0$-module up to nilpotents. Thus, the perfection of $0^*_{H_{\frakm_0}^d(R_0)}$ inside $H^d_{\frakm}(R)$ is a holonomic $R_0$-module whose image under $\text{Sol}$ is an \'etale sheaf supported on $\Spec(R/\frakm)$. Since $R/\frakm$ is algebraically closed, the condition on $\varphi$-fixed elements implies that $0^*_{H_{\frakm_0}^d(R_0)}$ is nilpotent.
\end{proof}

\begin{lemma}[Emerton--Kisin] \label{lemm--EKsimple}
    Let $X$ be a smooth irreducible scheme of finite type over a perfect field $k$. Let $\mathcal{L}$ be an \'etale local system of $\bbF_p$-vector spaces on $X$. Then $\mathcal{L}[\dim X]$ is simple as a perverse sheaf if and only if it is simple as a local system.
\end{lemma}

\begin{proof}
    The property of being perverse is \'etale-local so that Theorem \ref{thm--CM-conditions} implies $\mathcal{L}[\dim X]$ is perverse. The part about simplicity follows from the claim that every perverse subsheaf of $\mathcal{L}[\dim X]$ is again a shifted local system. Indeed, Gabber's result \cite[Corollary 12.4]{Gab04} implies every perverse subsheaf is constructible, and then there are multiple ways to proceed; here we sketch the argument of Emerton--Kisin in \cite[Corollary 4.3.3]{EK04b}. Via their Riemann--Hilbert correspondence, $\mathcal{L}[\dim X]$ corresponds to a unit $\varphi$-crystal, i.e., an $\mathcal{O}_{X}[\varphi]$-module $M$, locally free of finite rank over $\mathcal{O}_X$, where the unit condition means that the adjoint map $\varphi^*M \rightarrow M$ is an isomorphism. Their correspondence is a perverse t-exact anti-equivalence, so that perverse subsheaves of $\mathcal{L}[\dim X]$ correspond to unit $\mathcal{O}_{X}[\varphi]$-module quotients of $M$. Then the key input is that any such quotient is is locally free \cite[Proposition 1.2.3]{EK04b}, so that its Riemann--Hilbert partner is a shifted local system.
\end{proof}

Our main result in this subsection characterizes those schemes for which a simple local system corresponds to a simple perverse sheaf.

\begin{theorem} \label{prop-Fp-simple}
Let $k$ be a perfect field of characteristic $p$ and let $X$ be a connected scheme isomorphic to the perfection of a finite type $k$-scheme. Then the following are equivalent.
\begin{enumerate}
    \item $X$ is $\varphi$-rational in the sense of Definition \ref{def--rational}.
    \item The shifted constant sheaf $\bbF_p[\dim X]$ is a simple perverse sheaf.
\end{enumerate}
\end{theorem}

\begin{proof}
First suppose that $X$ is $\varphi$-rational. Then $\bbF_p[\dim X]$ is perverse by Theorem \ref{thm--CM-conditions}. Furthermore, $X$ is irreducible by Lemma \ref{lemm--rationalequiv} and \cite[Theorem 3.1]{DMP24}. Let $U \subset X$ be a nonempty open subscheme isomorphic to the perfection of a smooth finite type $k$-scheme, which exists by \cite[Tag 056V]{StaProj}. Let $i \colon X \setminus U \rightarrow X$ be a complementary closed immersion. Then $\bbF_p[\dim X]_U$ is perverse and simple on $U$ by Lemma \ref{lemm--EKsimple}, so it suffices to show the intermediate extension to $X$ is $\bbF_p[\dim X]$. Clearly $i^*\bbF_p[\dim X] \in {}^{p}D^{\leq -1}(X \setminus U, \mathbb{F}_p)$, and it remains to verify that $Ri^!\bbF_p[\dim X] \in {}^{p}D^{\geq 1}(X \setminus U, \mathbb{F}_p)$. Let $(R, \frakm)$ be a strict henselization of the local ring at a point in $X \setminus U$. By Lemma \ref{lemm--rationalequiv} the hypothesis of Lemma \ref{lemm--puncturerational} are satisfied, and in particular $R$ is a domain. We must therefore verify that $R^{i}\Gamma_{\{\frakm\}}(\bbF_p) = 0$ for $i \leq \dim R$, where $\bbF_p$ is viewed as an \'etale sheaf on $\Spec(R)$.  When $i < \dim R$ this vanishing follows from perversity, and furthermore $H^i_{\frakm}(R) = 0$ by Cohen--Macaulayness. The case $i = \dim R$ then follows from Proposition \ref{lemm--puncturerational}, by applying $R\Gamma_{\{\frakm\}}$ to the Artin--Schreier sequence \eqref{eq--AS}.

For the other direction, let $j \colon U \to X$ be an irreducible open subscheme isomorphic to the perfection of a smooth finite type $k$-scheme. By simplicity we must have $\bbF_p[\dim X] \cong j_{!*}(\bbF_p[\dim X]_U)$. On the other hand, $j_{!*}(\bbF_p[\dim X])$ is supported on the closure of $U$, so $X$ is irreducible. We now show by descending induction on $\dim \overline{\{ x\}}$ that the local ring at $x \in X$ is $\varphi$-rational. If $x \in U$ this follows since regular local rings $\varphi$-rational even before passing to the perfection \cite[Theorem 2.1 a)]{HH89}. If $x \in X \setminus U$ we may assume the punctured spectrum of the local ring at $x$ is $\varphi$-rational, so the same is true of the strict henselization $(R, \frakm)$ by Lemma \ref{lemm--rationalpuncture}.
Since $X$ is irreducible, the condition $j_{!*}(\bbF_p[\dim X]_U) =\bbF_p[\dim X]$ implies that $R^{i}\Gamma_{\{\frakm\}}(\bbF_p) = 0$ for $i \leq \dim R$, where $\bbF_p$ is viewed as an \'etale sheaf on $\Spec(R)$. Perversity of $\bbF_p[\dim X]$ implies that $H^i_{\frakm}(R) =0$ for $i < \dim R$ (Theorem \ref{thm--CM-conditions} and Lemma \ref{lemm--loccohstrict}). Then by the Artin--Schreier sequence \eqref{eq--AS}, $H^{\dim R}_{\frakm}(R)$ has trivial $\varphi$-invariants. The remaining hypotheses of Proposition \ref{lemm--puncturerational}, equidimensionality in particular, are satisfied for $R$ by Lemma \ref{lemm--invariance}. Thus, $R$ is $\varphi$-rational, and hence so is the local ring at $x$ by Lemma \ref{lemm--rationalequiv}.
\end{proof}

A priori, the fact that the irreducibility of a scheme is not an \'etale-local property could prevent the simplicity of $\bbF_p[\dim X]$ from being an \'etale-local property. However, the relation with $\varphi$-nilpotence shows that this is not the case.

 \begin{corollary}
     Let $k$ be a perfect field of characteristic $p$ and let $X_0$ be a finite type $k$-scheme. If the shifted constant sheaf $\bbF_p[\dim X_0]$ is a simple perverse sheaf, then $X_0$ is geometrically unibranch.
 \end{corollary}

 \begin{proof}
     By Theorem \ref{prop-Fp-simple} and Lemma \ref{lemm--rationalequiv}, the local rings of $X_0$ are $\varphi$-nilpotent, so the result follows from  \cite[Theorem 3.1, Remark 3.2]{DMP24}.
 \end{proof}

\begin{remark}
Let $X$ be a normal irreducible scheme of finite type over a perfect field $k$. Two important theorems in commutative algebra assert that the absolute integral closure $X^+$ of $X$ in its field of fractions is Cohen--Macaulay (due to Hochster--Huneke \cite{HH92}) and $\varphi$-rational (due to Smith \cite{Smi94}) in an appropriate sense. This is true more generally over a $\varphi$-finite base. As observed in \cite[\S 5.6]{BBL+23}, these theorems can be recovered from results such as Theorem \ref{prop-Fp-simple}. Informally, the idea is to show that $\bbF_p[\dim X]_{X^+}$ is a simple perverse sheaf on $X^+$. To prove this, one must compute $i^*\bbF_p[\dim X]_{X^+}$ and $Ri^! \bbF_p[\dim X]_{X^+}$, where $i \colon Z \rightarrow X^+$ is a proper closed subscheme. But the $*$-pullback is constant, and the $!$-pullback vanishes since $\bbF_p[\dim X]_{X^+}$ is the $*$-extension of its restriction to any open subset  \cite[Proposition 3.10]{Bha20}. We refer to  \cite[\S 5.6]{BBL+23} for more details. 
\end{remark}

\section{Global +-regularity and inversion of adjunction}

In this section, we review some of the material from \cite{BMP+23} on globally +-regular varieties, explain the proof of their criterion for inversion of adjunction, and adapt it to a certain asymptotic analog. This is going to be applied later to certain Demazure varieties.

\begin{remark}
	In positive characteristic, ideas such as these have been known for at least a decade in advance. For instance, Das \cite{Das15} proved inversion of adjunction for strong $\varphi$-regularity in characteristic $p$, but we would rather avoid his treatment, because it circumvents the Kawamata--Viehweg $+$-vanishing of \cite{Bha12} and because it forces us to work with $\mathbb{Z}_{(p)}$-divisors everywhere instead of $\mathbb{Q}$-divisors. 
\end{remark} 

\subsection{Global $+$-regularity}
Let $k$ be a perfect field of characteristic $p$ and $X$ be a finite type connected normal $k$-scheme. It is helpful to consider the notion of boundary and subboundary $\mathbb{Q}$-divisors, as they constitute a very flexible tool in studying singularities. 

\begin{definition}
	Let $\Delta=\sum_i r_i D_i$ be a $\bbQ$-divisor on $X$, i.e., a finite rational linear combination of prime divisors on $X$. We say that $\Delta$ is a boundary (resp.~a subboundary) if $0\leq r_i \leq 1$ for all $i$ (resp.~ if $0\leq r_i <1$). We refer to $(X,\Delta)$ as a boundary (resp.~subboundary) pair.
\end{definition}

Let us recall the notion of global $+$-regularity following \cite[Definition 6.1]{BMP+23}.

\begin{definition} \label{def:G+R}
	We say that the pair $(X,\Delta)$ is globally $+$-regular if for every finite cover $f\colon Y \to X$ with $Y$ connected normal, the natural map $\calO_X \to f_*\calO_Y( \lfloor f^*\Delta \rfloor )$ splits in the category of $\calO_X$-modules.
\end{definition}

If $\Delta=0$, then we simply say that $X$ is globally $+$-regular. Note that this condition only has to be verified for a cofinal family of finite covers $f$. 
By the cyclic covering trick, we may even assume that $f^*\Delta$ is integral, compare with \cite[Remark 6.2]{BMP+23}.
Let us start with the first basic stability property.

\begin{lemma}
	If the boundary pair $(X,\Delta)$ is globally $+$-regular, the same holds true for $(X,\Delta')$ for any boundary $\Delta'\leq \Delta$.
\end{lemma}

\begin{proof}
	Compose the inclusion $f_*\calO_Y(\lfloor f^*\Delta' \rfloor)  \to f_*\calO_Y(\lfloor f^*\Delta \rfloor) $ with the given section of $\calO_X \to f_*\calO_Y(\lfloor f^*\Delta \rfloor)$.
\end{proof}

In particular, $(X,\Delta)$ being globally $+$-regular implies that $(X,\epsilon\Delta)$ also is for every $0\leq \epsilon \leq 1$. More importantly, global $+$-regularity satisfies proper descent:

\begin{proposition}
	Let $f \colon X \to Y$ be a proper birational map of normal connected $k$-schemes of finite type. If $(X,\Delta)$ is globally $+$-regular, then so is $(Y,f_*\Delta)$.
\end{proposition}

\begin{proof}
	This is \cite[Proposition 6.19]{BMP+23}. Note that pushforwards and pullbacks of $\mathbb{Q}$-divisors along alterations of normal connected finite type are defined locally in codimension $1$ on principal divisors via the norm map and the inclusion map, respectively, see \cite[Tag 02RS]{StaProj}, and then extended by normality to the entire space. Let $g\colon Z \to Y$ be a finite cover by a normal integral $k$-scheme such that $g^*f_*\Delta$ is an integral divisor. Let $W$ be the normalization of $X\times_Y Z$ with base maps $g'\colon W \to X$ and $f'\colon W\to Z$. Then, we know that $\calO_X \to g'_*\calO_W(g'^*\Delta)$ splits in $\calO_X$-modules. The same holds therefore for $\calO_Y \to g_*\calO_Z( f'_*g'^*\Delta )$ in $\calO_Y$-modules, because there is a natural map $g_*\calO_Z(f'_*g'^*\Delta)\to f'_*g'_*\calO_W(g'^*\Delta)$. Noticing that $g^*f_*\Delta= f'_*g'^*\Delta$, we deduce our desired splitting.
\end{proof}

Next, we translate the notion of globally $+$-regularity in terms of trace maps by applying Grothendieck--Serre duality. Recall that there is a 6-functor formalism on the category of quasi-coherent sheaves and the canonical sheaves $\omega_X$ arise as the $H^{-\mathrm{dim}(X)}$ of the complex $Rp^!k$, where $p\colon X\to \mathrm{Spec}(k)$ denotes the structure morphism. Since we assume $X$ to be normal and connected, it turns out that $\omega_X$ is a reflexive sheaf and we usually fix an arbitrary canonical divisor $K_X$ such that $\omega_X\simeq \mathcal{O}_X(K_X)$. 

\begin{proposition} \label{prop--Trace}
	Assume $K_X+\Delta$ is $\mathbb{Q}$-Cartier. Then, the boundary pair $(X,\Delta)$ is globally $+$-regular if and only if the trace map \begin{equation} \label{eq--Trace}
		H^0(Y,\calO_Y(K_Y-\lfloor f^*(K_X+\Delta)\rfloor)) \to H^0(X,\calO_X)
	\end{equation}
is surjective for all connected normal finite covers $f\colon Y \to X$.
\end{proposition}

\begin{proof}
	This is a particular case of \cite[Proposition 6.8]{BMP+23}. By duality, the natural map $\calO_X \to f_*\calO_Y(\lfloor f^*\Delta\rfloor)$ of $\calO_X$-modules is a split injection if and only if the map
	\begin{equation}
		\calH\mathrm{om}_{\calO_X}(f_*\calO_Y(\lfloor f^*\Delta\rfloor),\calO_X)\to \calO_X
	\end{equation}
of $\calO_X$-modules is a split surjection. By Grothendieck--Serre duality, the left side identifies with $f_*\mathcal{H}\mathrm{om}_{\mathcal{O}_Y}(\mathcal{O}_Y(\lfloor f^*\Delta \rfloor) ,f^!\calO_{X})$ where $f^!$ is the abelian truncation of the shriek pullback $Rf^!$. Note that this is reflexive, so to compute it we are allowed to restrict to the smooth locus of $Y$. Over there, $\lfloor f^*\Delta \rfloor$ becomes an actual Cartier divisor, so, in particular, we can write the left side as $f_*f^!\calO_X(-\lfloor f^*\Delta \rfloor)$ by pulling the divisor across the Hom. On the other hand, by definition of the canonical divisor, we have $f^!\calO_X(K_X)=\calO_Y(K_Y)$, and since we are over the smooth locus of $Y$, we get the identity $f^!\mathcal{O}_X=\calO_Y(K_Y-K_X)$ and our sheaf identifies with $f_*\calO_Y(K_Y-\lfloor f^*(K_X+\Delta)\rfloor)$, just like in the statement of the proposition. Now, since $\calO_X$ is free, the surjectivity of the trace map can be tested at the level of global sections.
\end{proof}

Motivated by the previous proposition, one has the $k$-module of $+$-stable sections $B^0(X,\Delta;\calO_X)$ in \cite[Definition 4.2]{BMP+23} given by the intersection across all normal finite covers $f\colon Y \to X$ of the images of the trace maps \eqref{eq--Trace} appearing in the statement of Proposition \ref{prop--Trace}. In particular, global $+$-regularity amounts to demanding an equality $B^0(X,\Delta,\calO_X)=H^0(X,\calO_X)$. We finish this subsection with the following quite non-standard notion:

\begin{definition}
	We say that a boundary pair $(X,\Delta)$ is $\mathbb{Q}$-Fano if the $\mathbb{Q}$-divisor $K_X+\Delta$ is $\mathbb{Q}$-Cartier and anti-ample.
\end{definition}

Eventually, we will want to provide an inductive criterion for lifting global $+$-regularity along closed subschemes and this positivity condition will play a significant role.

\subsection{Pure variant and inversion of adjunction}
Our next topic consists of a variant of global $+$-regularity defined along a prime divisor $S \subset X$. The following notion is a simplification of \cite[Definition 6.24]{BMP+23}. Let $(X,\Delta)$ be a pair consisting of a finite type $k$-scheme and an effective $\bbQ$-divisor $\Delta$ on $X$. Assume $\Delta=S+B$ where $S$ is a prime divisor and $B$ an effective $\bbQ$-divisor on $X$ with irreducible components different from $S$.

\begin{definition}
	 We say that $(X,S+B)$ is purely globally $+$-regular along $S$ if the map of $\calO_X$-modules $\calO_X \to f_*\calO_Y(-S_Y+\lfloor f^*(S+B)\rfloor)$ splits for every finite cover $f\colon Y \to X$ with $Y$ connected normal, where the $S_Y \subset Y$ form a compatible family of prime divisors lying over $S\subset X$.
\end{definition}

Using the Galois action on $X^+$ over $X$, we can show that the previous definition is independent of the choice of the prime divisors $S_Y \subset Y$ (equivalently, of an absolute integral closure $S^+\subset X^+$). It also amounts to purity of the map $\mathcal{O}_X \to \calO_X^+(-S^++\pi^*(S+B))$, where $\pi \colon X^+\to X$ is the absolute integral closure. There is a close relationship between pure global $+$-regularity and global $+$-regularity after slightly tweaking the divisors.

\begin{lemma}\label{lem_pure_to_non_pure_plus_regularity}
	If $(X,S+B)$ is purely globally $+$-regular along $S$, then $(X,(1-\epsilon) S+B)$ is globally $+$-regular for every rational number $0< \epsilon \leq 1$.
\end{lemma}
\begin{proof}
	This is \cite[Lemma 4.26]{BMP+23}. We just have to notice that $f^*(\epsilon S+B)\leq -S_Y+f^*(S+B)$ for sufficiently large normal finite covers $f\colon Y\to X$, so that the pure $+$-splitting of $(X,S+B)$ along $S$ factors over a $+$-splitting for the pair
$(X,(1-\epsilon) S+B)$.\end{proof}

\begin{proposition}
	The pair $(X,S+B)$ is purely globally $+$-regular along $S$ if and only if the trace map \begin{equation}
		H^0(Y,\calO_Y(K_Y+S_Y-\lfloor f^*(K_X+S+B)\rfloor)) \to H^0(X,\calO_X)
	\end{equation}
	is surjective for all normal finite covers $f\colon Y\to X$.
\end{proposition}

\begin{proof}
	The proof is the same as the non-pure along $S$ case, requiring us to check that the $\calO_X$-module dual of $f_*\calO_Y(-S_Y+\lfloor f^*(S+B) \rfloor)$ equals $f_*\calO_Y(K_Y+S_Y-\lfloor f^*(K_X+S+B)\rfloor)$.
\end{proof}

Again, there exists a module $B^0_S(X,S+B;\calO_X)$ of pure $+$-stable sections along $S$, see \cite[Definition 4.21]{BMP+23}, and pure global $+$-regularity along $S$ translates into an equality $B^0_S(X,S+B;\calO_X)=H^0(X,\calO_X)$. 
The pure along $S$ variant of global $+$-regularity was set up in this way, precisely because we want to study how to lift global $+$-regularity from a prime divisor $S$ to the whole $k$-variety $X$ -- this is known as inversion of adjunction.

\begin{theorem}[\cite{BMP+23}]\label{thm_inversion_adjunction}
	Let $X$ be a connected normal proper $k$-scheme, $S\subset X$ a normal prime divisor, and $B$ a subboundary with components different from $S$. If $(X,S+B)$ is $\mathbb{Q}$-Fano, then $(X,S+B)$ is purely globally $+$-regular along $S$ if and only if $(S,B|_S)$ is globally $+$-regular.
\end{theorem}

\begin{proof}
This is a particular case of \cite[Theorem 7.2]{BMP+23}, see also \cite[Corollary 7.5]{BMP+23}, and we give a sketch of the argument. During the proof, we use the shorthand $\Delta=S+B$. By Serre duality, we can identify the trace maps with the natural maps
	\begin{equation}
		H^d(X,\calO_X(K_X))\to H^d(Y,\calO_Y(\lfloor f^*(K_X+\Delta)\rfloor))
	\end{equation}
induced by pullback along $f^*$ and multiplication by the divisor $\lfloor f^*\Delta \rfloor$, and similarly
	\begin{equation}
	H^d(X,\calO_X(K_X))\to H^d(Y,\calO_Y(-S_Y+\lfloor f^*(K_X+\Delta)\rfloor))
\end{equation}
in the pure along $S$ case. Note that $\calO_X(K_X+S)$ pulls back to $\calO_S(K_S)$ because this holds away from codimension $2$, and then we apply Hartogs' theorem by normality of $S$ and $X$. In particular, the associated long exact sequence yields a connecting homomorphism
\begin{equation}
	H^{d-1}(S,\calO_S(K_S))\to H^d(X,\calO_X(K_X))
\end{equation}
which is surjective by normality and connectedness: indeed, it arises by dualising on $k$-modules the non-zero ring map $H^0(X,\mathcal{O}_X)\to H^0(S,\mathcal{O}_S)$ between finite field extensions of $k$.
Similarly, we can connect the right sides of the pullback maps via the following map
\begin{equation}
	H^{d-1}(S^+,\calO_S^+(\nu_S^*(K_S+B|_S))) \to H^d(X^+,\calO_X^+(-S^++ \nu_X^*(K_X+\Delta)))
\end{equation}
where we let the $+$-notation denote the colimit with respect to a family of connected normal finite covers $f\colon Y\to X$, and $\nu$ is the structure map of the absolute integral closures. The kernel of the connecting homomorphism at the $+$-level is given by the image of $H^{d-1}(X^+,\calO_X^+(\nu_X^*(K_X+\Delta)))$. The latter vanishes by anti-ampleness of the $\mathbb{Q}$-Cartier divisor $K_X+\Delta$ and the Kodaira $+$-vanishing theorem of \cite{Bha12}. A diagram chase reveals that injectivity for $S$ gives rise to injectivity for $X$, and vice-versa.
\end{proof}

\subsection{An asymptotic variant}

\Cref{thm_inversion_adjunction} provides a criterion for inversion of adjunction of global $+$-regularity but has the somewhat unpleasant feature that it lifts global $+$-regularity to at most pure global $+$-regularity. At the same time, the latter comes pretty close to global $+$-regularity itself by \Cref{lem_pure_to_non_pure_plus_regularity}. This leads us to formulate a variant that treats boundary pairs asymptotically and improves the clarity of our exposition when applying the criterion to Demazure varieties.

\begin{definition}
Given a boundary decomposition $\Delta=S+B$ with $S$ prime and $B\geq 0$ with no common components with $S$, we similarly say that the boundary pair $(X,\Delta)$ is asymptotically purely $\mathbb{Q}$-Fano along $S$ if there exist arbitrarily close subboundaries $B'< B$ such that $K_X+S+B'$ is an anti-ample $\mathbb{Q}$-Cartier $\mathbb{Q}$-divisor. 
\end{definition} 

The definition above is again quite non-standard, but it fits well within our paper. The next step is to define the asymptotic analog of global $+$-regularity.

\begin{definition}
	We say that the boundary pair $(X,\Delta)$ is asymptotically globally $+$-regular if for all subboundaries $\Delta'<\Delta$, the pair $(X,\Delta')$ is globally $+$-regular in the usual sense. 
\end{definition}

Note that here the condition applies to all smaller subboundaries, because global $+$-regularity is stable under parallelipipeds, unlike ampleness. We can safely ignore a corresponding asymptotic notion of pure global $+$-regularity along a prime divisor, as the criterion for inversion of adjunction now takes the following form.

\begin{corollary} \label{cor:asymptotic}
	Let $X$ be a connected normal proper $k$-scheme, $S\subset X$ a normal prime divisor, and $B$ a boundary with components different from $S$. If $(X,S+B)$ is asymptotically purely $\mathbb{Q}$-Fano along $S$ and $(S,B|_S)$ is asymptotically globally $+$-regular, then $(X,S+B)$ is asymptotically globally $+$-regular.
\end{corollary}

\begin{proof}
	Let $B'<B$ be a subboundary such that the corresponding pair $(X,\Delta')$ with $\Delta'=S+B'$ is $\mathbb{Q}$-Fano. Now, since we know that $(S,B'|_S)$ is globally $+$-regular, we may apply \Cref{thm_inversion_adjunction} to get that $(X,\Delta')$ is purely globally $+$-regular along $S$. But then $(X,(1-\epsilon)S+B')$ is actually globally $+$-regular for any $\epsilon>0$. Letting $\epsilon$ go to $0$ and $B'$ to $B$, we get arbitrarily close to the original boundary $\Delta$, so it induces an asymptotically globally $+$-regular pair.
\end{proof}

\begin{remark}
	There is a corresponding version of the corollary which is an equivalence between the behavior of the pairs $(X,\Delta)$ and $(S,B|_S)$, but it requires defining asymptotic pure global $+$-regularity along a divisor. The only thing happening here is that the asymptotic pure version for $(X,\Delta)$ implies the asymptotic non-pure one for the same pair without tampering with the divisor, precisely because of the asymptoticity.
\end{remark}

\section{Abstract Demazure and Schubert varieties}\label{sec_abstract_BSDH}
Several decades ago, Bott--Samelson, Demazure and Hansen independently discovered certain desingularizations of closed equivariant subvarieties of finite flag varieties. Since their first appearance, many other analogs have been constructed using different incarnations of group theory, such as Kac--Moody algebras, loop groups, etc.
In this section, we discuss an axiomatic setup for these types of varieties, which we call Bott--Samelson--Demazure--Hansen varieties (BSDH varieties for short), that encapsulates every single example we know of appearing in nature, and simultaneously does not rely on any group theory.

\subsection{Basic properties} We start with the key definition of this section and study its most immediate properties.

\begin{definition}
A BSDH variety over a field $k$ is a variety $X$ endowed with an inductive structure depending on its dimension as follows. For $\dim(X) = 0$, we require $X = \operatorname{Spec}(k)$. For $\dim(X) > 0$, this structure is given by a Zariski locally trivial $\mathbb{P}^1_k$-fibration $p \colon X \to Y$ equipped with a distinguished section, such that the base $Y$ is also a BSDH variety.
\end{definition}

\begin{remark}
    It is well-known that a Zariski locally trivial $\mathbb{P}_k^1$-fibration over a regular noetherian base is the projectivization $\mathbb{P}(\mathcal{E})$ of a rank $2$ vector bundle $\mathcal{E}$, e.g.~\cite[Ex.~7.10(c)]{Har77}. A section of $\mathbb{P}(\mathcal{E})$ is the same as a rank $1$ quotient of $\mathcal{E}$. Thus, a BSDH $k$-variety can also be described as a finite tower of varieties $X_k \rightarrow X_{k-1} \rightarrow \cdots \to X_0=\Spec(k)$, where each map is the projectivization of a rank $2$ vector bundle extension $\mathcal{L} \rightarrow \mathcal{E} \rightarrow \mathcal{O}$. The case where each rank $2$ bundle splits as a direct sum is called a Bott tower, see e.g.~\cite{GK94} for the complex-analytic case and the references therein. The $\bbP^1$-fibration $p\colon X\to Y$ for a Bott tower admits instead two disjoint sections.
\end{remark}

For a BSDH $k$-variety $X$ of dimension $n$, we can inductively construct closed immersions $\iota_J \colon X_J \to X$ for any subset $J \subset [n]:=\{1,\dots,n\}$ and projections $p_m\colon X\to X_m:=X_{[m]}$ for all $m\leq n$. We call $X_{\emptyset} = \Spec(k)$ the origin.
We have obvious simple normal crossings divisors $X_{[n]\setminus i}\subset X$ denoted by $D_i$. We denote the curves $X_{\{i\}}\subset X$ by $C_i$, and we note that there are isomorphisms $C_i \simeq \mathbb{P}^1_k$. 
The Picard group of $X$ can then be computed as follows.

\begin{proposition}
Let $X$ be a BSDH $k$-variety. Then there is an isomorphism $\mathrm{deg} \colon \mathrm{Pic}(X) \xrightarrow{\sim} \mathbb{Z}^{\oplus n}$ given by the degree of the restriction to $C_i$ for every $i=1,\dots,n$.
\end{proposition}

\begin{proof}
	We simply follow the proof of \cite[Lemma 4.7]{FHLR25}: the Zariski fibration gives us a short exact sequence $0\to \mathrm{Pic}(Y)\to \mathrm{Pic}(X)\to \mathrm{Pic}(\mathbb{P}^1_k)\to 0$ that splits due to the existence of a section $s\colon Y\to X$. Retracing the inductive construction of the $C_i$, this yields the claim.
\end{proof}

From now on, we will denote the preimage of $(d_1,\dots,d_n)\in \mathbb{Z}^n$ under the above isomorphism $\mathrm{deg}$ by $\mathcal{O}_X(d_1,\dots,d_n)$. Given a $\bbQ$-divisor $D \in \mathrm{Div}(X)$, we define its degree $\mathrm{deg}(D)\in \bbZ^n$ as the only sequence such that $\calO_X(mD)\simeq \calO_X(m\deg(D))$ for all sufficiently divisible $m$. Let $\theta_X$ be a $\bbQ$-divisor such that $\deg(\theta_X)=(1,\dots,1)$ and denote $\partial X=\sum_{i=1}^n D_i$. We can compute the canonical divisor of a BSDH $k$-variety as follows.

\begin{proposition}
Let $X$ be a BSDH $k$-variety and $\theta_X$ be a $\bbQ$-divisor such that $\deg(\theta_X)=(1,\dots,1)$. Then, the divisor $\partial X +\theta_X$ is an anti-canonical divisor for $X$.	
\end{proposition}

\begin{proof}
	We have to calculate $\omega_X$. We have an adjunction formula for canonical divisors along regular immersions, so that $\omega_{X}(\partial X)=\omega_{D_i}(\partial {D_i})$. Continuing inductively, we see that the restriction of $\omega_X(\partial X)$ to $C_j$ equals $\mathcal{O}_{C_j}(-1)$. In particular, it follows that $\omega_X\simeq \calO_X(-\partial X-\theta_X)$, as desired.
\end{proof}

Before moving on, we need the following lemma on $\mathbb{Q}_{>0}$-semigroups which will help us to  produce ample perturbations without explicitly calculating the ample cone. We first introduce the following notation. 
For any square matrix $A$, let $A \cdot \mathbb{Q}_{>0}^n \subset \mathbb{Q}^n$ be the $\mathbb{Q}_{>0}$-semigroup generated by its columns.

\begin{lemma}\label{lem_intersecting_monoids}
    Let $A, B \in \mathrm{GL}_n(\mathbb{Q})$ be two upper triangular matrices with positive diagonal entries. Then the $\bbQ_{>0}$-semigroups $A \cdot \mathbb{Q}_{>0}^n$ and $B \cdot \mathbb{Q}_{>0}^n$ have non-empty intersection.
\end{lemma}

\begin{proof} We induct on the $n$, with the case $n=1$ being trivial. If $n > 1$, let $A'$ and $B'$ be the minors obtained by removing the first column and row of $A$ and $B$, respectively. By induction, there exist $\mu_i > 0, \lambda_i > 0$ for $i=1, \ldots, n-1$ such that $\sum_{i=1}^{n-1} \mu_i \text{col}_i(A') = \sum_{i=1}^{n-1} \lambda_i \text{col}_i(B')$. Then \begin{equation}\sum_{i=2}^{n} \mu_{i-1}\text{col}_i(A) - \sum_{i=2}^{n} \lambda_{i-1}\text{col}_i(B) = \begin{pmatrix} r \\ 0 \\ \vdots \\ 0\end{pmatrix}\end{equation} for some $r \in \mathbb{Q}$. Let $a, b \in \mathbb{Q}_{>0}$ be the top left entries of $A$ and $B$, respectively. Pick $\mu_0, \lambda_0 \in \mathbb{Q}_{>0}$ such that $\mu_0 a - \lambda_0 b = r$. Then $\sum_{i=1}^{n} \mu_{i-1} \text{col}_i(A) = \sum_{i=1}^{n} \lambda_{i-1} \text{col}_i(B)$.
    \end{proof}

Usually, when working with classical Demazure varieties, we can identify the ample cone as given by strictly dominant vectors (i.e., strictly decreasing sequences). This is not necessarily true in our setup and we can avoid computing the ample cone thanks to the previous lemma. Let us now illustrate its usefulness.

\begin{lemma}\label{lem_intersection_divisor_ample_cones}
	Let $X$ be a BSDH $k$-variety equipped with its distinguished divisors $D_i$ and a $\bbQ$-divisor $\theta_X$ such that $\mathrm{deg}(\theta_X)=(1,\ldots,1)$. Then the following subsets meet the ample cone:
    \begin{enumerate}
        \item The $\bbQ_{>0}$-semigroup generated by $D_i$ for $1\leq i \leq n$.
        \item The $\bbQ_{>0}$-semigroup generated by $D_i$ for $1\leq i <n$ and $\theta_X$.
    \end{enumerate}
\end{lemma}

\begin{proof}
We deduce this from Lemma \ref{lem_intersecting_monoids}. First, we consider the pullback $A_i$ along $X\to X_i$ of an ample divisor bundle on $X_i$ (the former exists because the latter is a projective $k$-variety). Next, we look at the associated square matrix $A$ with columns $\mathrm{deg}(A_1),\dots, \mathrm{deg}(A_n)$. Note that the preimage of the semigroup $A\cdot \bbQ_{>0}^n$ along the degree isomorphism $\mathrm{deg}\colon \mathrm{Pic}(X)\to \bbZ^n$ is entirely contained in the ample cone, by definition. Also by construction, the pullback of $A_i$ to $C_j$ for $j>i$ is trivial, as $C_j\to X_i$ is the constant map to the origin, so $A$ is upper triangular. Similarly, the pullback to $C_i$ has positive degree, so its diagonal has strictly positive entries. Thus, this matrix fits into the framework of our previous lemma.

Next, we consider the square matrices $B$ and $C$ whose first $n-1$ columns are given by $\mathrm{deg}(D_1),\dots,\mathrm{deg}(D_{n-1})$, whereas the last one is respectively $\mathrm{deg}(D_n)$ or $\mathrm{deg}(\theta_X)$. Clearly, the divisors $D_i$ are obtained by pullback along $X\to X_i$, so we deduce that $B$ and $C$ are upper triangular. Finally, both $D_i$ and $\theta_X$ have degree $1$ upon restricting to $C_i$, so the diagonals of $B$ and $C$ have strictly positive entries. Now, we notice that the semigroups $B\cdot \bbQ^n_{>0}$ and $C\cdot \bbQ^n_{>0}$ are the image under the degree map of the semigroups in the statement of the lemma. Thanks to Lemma \ref{lem_intersecting_monoids}, we now know that $A\cdot \bbQ^n_{>0}$ meets both $B\cdot \bbQ_{>0}^n$ and $C\cdot \bbQ_{>0}^n$. The statement of the lemma now follows. 
\end{proof}

\subsection{Global $+$-regularity of based BSDH $k$-varieties}

In this subsection, we introduce the missing hypothesis necessary to tame the singularities of BSDH $k$-varieties, that relates to the behavior of theta divisors. We refer to these varieties as based BSDH $k$-varieties and then we will prove that they are globally $+$-regular.

\begin{definition}
	A based BSDH $k$-variety $X$ is a BSDH $k$-variety such that the origin does not lie in the stable base locus of the line bundle $\mathcal{O}_X(1,\dots,1)$.
\end{definition}

Note that we do not require $\mathcal{O}_X(1,\dots,1)$ to be semi-ample or big. We already know that the canonical divisor of $X$ is of the form $\partial X+\theta_X$, but for a BSDH $k$-variety we may and do assume that the theta divisor $\theta_X$ does not vanish at the origin. Indeed, the condition on the stable base locus is equivalent to finding an origin-avoiding effective theta $\bbQ$-divisor $\theta_X$ on $X$. 

\begin{remark}
    Bott towers $X$ are always based with a possible choice of $\theta_X$ being given by the sum of the north pole sections, see \cite[Lemma 3.5]{NC18}. In particular, our Theorem \ref{thm_globally_+_regular_abstract_demazure} below recovers the known result by Smith that projective toric $k$-varieties are globally $\varphi$-regular, see \cite[Proposition 6.4]{Smi00}. Before this paper, every group-theoretic BSDH $k$-variety which we are aware of was based with a possible choice of $\theta_X$ arising from taking global sections of a critical line bundle on the corresponding flag variety. Later we will encounter BSDH $k$-varieties whose stable base locus of $\calO_X(1,\ldots,1)$ does not avoid the origin.
\end{remark}

We fix a choice of origin-avoiding effective $\bbQ$-divisor $\theta_X$ once and for all
and we now prove that BSDH $k$-varieties are globally $+$-regular.

\begin{theorem}\label{thm_globally_+_regular_abstract_demazure}
Let $X$ be a based BSDH $k$-variety and $\theta_X$ be an origin-avoiding effective $\bbQ$-divisor on $X$. The pair $(X,\partial X+\theta_X)$ is asymptotically globally $+$-regular.
\end{theorem}

\begin{proof}
We set $S:=D_n=X_{n-1}$ and $\Delta_X:=\partial X +\theta_X$ and our goal is to show that $(X,\Delta_X)$ is asymptotically purely Fano along $S:=D_n$, so that we may apply Corollary \ref{cor:asymptotic}, which is the asymptotic analog of \Cref{thm_inversion_adjunction}. The fact that $\theta_X$ is origin-avoiding enters the picture to ensure that its restriction to $S$ still differs from the restriction of $D_i$ for all $i <n$. 
	Thus, we have to slightly perturb the coefficients of $\Delta_X=\theta_X+\sum_{i\leq n}D_{i}$ to get a $\bbQ$-divisor 
	\begin{equation}\Delta_X'=r_0\theta_X+\sum_{i\leq n} r_{i}D_{i}\end{equation} with $r_i$ smaller but arbitrarily close to $1$ and $r_n=1$ so that $K_X+\Delta_X'$ is $\mathbb{Q}$-Cartier (trivial as $X$ is smooth) and anti-ample.
	For convenience, we set $\epsilon_i:=1-r_i$. Then, we deduce that
	\begin{equation}
		-K_X-\Delta_X'=\epsilon_0 \theta_X+\sum_{i}\epsilon_iD_{i}
	\end{equation}
	and it is enough to choose the $\epsilon_i$ such that the associated divisor is ample. For this to be possible, the $\mathbb{Q}_{> 0}$-semigroup spanned by $\theta_X$ and the $D_i$ for $i<n$ must meet the ample cone. But we verified this in the second half of Lemma \ref{lem_intersection_divisor_ample_cones}: 
    since the intersection of these two cones is stable under $\mathbb{Q}_{>0}$-homothety, we can choose the $\epsilon_i$ arbitrarily close to $0$ and such that the right side in the equation above is ample. Thus, we get the required asymptotical pure Fano property. 
\end{proof}

Furthermore, we can even deduce the more general property of being globally $\varphi$-regular from our proof. Recall that $X$ is globally $\varphi$-regular in the sense of \cite[Definition 3.5]{Smi00} if for any effective $\bbZ$-divisor $D$, there exists $e\gg 0$ such that $\calO_X\to \varphi^e_\ast\calO_X(D)$ splits in $\calO_X$-modules. More generally, a divisor pair $(X,\Delta)$ is globally $\varphi$-regular in the sense of \cite[Definition 3.1]{SS10} if instead $\calO_X\to \varphi^e_\ast\calO_X(\lceil (p^e-1)\Delta \rceil+D)$ splits in $\calO_X$-modules for $e\gg 0$. We will say that the divisor pair $(X,\Delta)$ is asymptotically globally $\varphi$-regular if it can be approximated by globally $\varphi$-regular pairs.

\begin{corollary} \label{cor--splitBSDH}
	Let $X$ be a based BSDH $k$-variety. Then, $(X,\partial X+\theta_X)$ is asymptotically globally $\varphi$-regular and compatibly $\varphi$-split with the BSDH subvarieties $X_J$ for all $J\subset [n]$.
\end{corollary}

\begin{proof}
	Global $+$-regularity applied to $(1-p^{-1})D_i$ implies that the map of coherent sheaves $\mathcal{O}_{X}\to \varphi_*\mathcal{O}_{X}((p-1)D_i)$
	splits in $\mathcal{O}_X$-modules. Twisting it by $\mathcal{O}(-D_i)$ and applying the projection formula, we deduce that $X$ is compatibly $\varphi$-split with each irreducible component $D_i$ of its boundary $\partial X$, compare with \cite[Theorem 1.4.10]{BK07}. This proves by induction that it is compatibly $\varphi$-split with any intersection and unions of those. 
	
	Next, we handle asymptotic global $\varphi$-regularity. We fix the divisor $\Delta=(\partial X+\theta_X)$ and approximate it by the convergent sequence $\frac{r-2}{r-1}\Delta$, where $r=p^c$ for some $c\in \bbZ_{>0}$ and $c$ goes to $\infty$. Let $A$ be an ample effective $\mathbb
{Z}$-divisor on $X$ lying in the $\mathbb{Q}_{>0}$-semigroup spanned by the $D_i$ for $i=1,\dots n$. Again, we know that such an $A$ exists by the first half of Lemma \ref{lem_intersection_divisor_ample_cones}. Its complement is an affine space, so it is locally $\varphi$-regular and we choose $q=p^e\gg 0$ with $c\mid e$ so that $r-1 \mid q-1$ and $q^{-1}A+\frac{r-2}{r-1}\Delta$ is a subboundary, i.e., less than or equal to $\Delta$. Then, we can guarantee the existence of a splitting \begin{equation}\mathcal{O}_{X}\to \varphi^e_\ast\mathcal{O}_X(A+\frac{(q-1)(r-2)}{r-1}\Delta) \end{equation}
	by global $+$-regularity and the inequality $\frac{(q-1)(r-2)}{r-1}\leq \lfloor \frac{q(r-2)}{r-1}\rfloor$. Now we can apply \cite[Theorem 3.9]{SS10} 
    to deduce that $X$ is globally $\varphi$-regular.
\end{proof}

\begin{remark}
   There is also a notion of pure global $\varphi$-regularity along $S$ for pairs $(X,S+B)$ and we can prove similar results for based BSDH $k$-varieties following Theorem \ref{thm_inversion_adjunction}, but we do not pursue this here.
\end{remark}

Let $\mathcal{L}$ be a semi-ample line bundle on $X$ and denote by $S$ its Stein factorization. This can be defined by first considering the natural map $\phi_n\colon X\to \bbP(H^0(X,\calL^{\otimes n}))$ for a sufficiently divisible $n$ such that $\calL^{\otimes n}$ becomes basepoint-free. Then, we take $S$ to be the Stein factorization of $\phi_n$. This coincides with the scheme-theoretic image of $\phi_n$ for all sufficiently divisible $n$ by \cite[Theorem 2.1.27]{Laz06}, so it is independent of $n$. Another way of constructing $S$, following \cite[Proposition 3.2]{AGLR22}, is to define $|S|$ as the quotient of $|X|$ given by collapsing connected $k$-subschemes over which $\calL$ is torsion, and taking $\calO_S$ as the (non-derived) pushforward of $\calO_X$ along $|X|\to |S|$.

Now, consider the closed embeddings $X_J\subset X$. The restriction of $\mathcal{L}$ to $X_J$ gives rise to a Stein factorization $X_J\to S_J$. There are natural rational maps $\bbP(H^0(X_J,\calL^n))\dashrightarrow \bbP(H^0(X,\calL^n))$, i.e., defined over an open subset. Since $X_J$ lands in the definition locus, we get natural maps $S_J\to S$ of Stein factorizations for the corresponding pullback of $\mathcal{L}$ to $X_J\subset X$. 

\begin{proposition} \label{Fregular-prop_abs}
The normal $k$-variety $S$ is globally $\varphi$-regular and the maps $S_J\to S$ are compatibly $\varphi$-split for all $J\subset [n]$. If $S_J\to S$ is a universal homeomorphism onto the image, then it is also a closed immersion. If $\mathcal{L}$ is big, then $X\to S$ is a rational resolution.

\end{proposition}
\begin{proof}
	Global $\varphi$-regularity is preserved along proper maps by \cite[Lemma 1.2]{LRPT06}, so the first claim follows from Corollary \ref{cor--splitBSDH}. For the second claim, it is crucial to show that the induced map $\mathcal{O}_{S}\to \mathcal{O}_{S_J}$ is surjective. This is equivalent to showing that the $\mathcal{O}_X(-D_i)$ has vanishing higher direct images along $X\to S$. By assumption, the proper map $S_J\to S$ is a universal homeomorphism onto the image, so it becomes a closed immersion after passing to perfections, yielding vanishing of the higher direct images at the perfect level.
    We can therefore descend it to our deperfection by using a $\varphi$-splitting of $X$ compatible with $D_i$. 
	
	As for rationality of $X\to S$, we can prove it in the same manner: we know that the higher direct images of the structure sheaf vanish at the level of absolute integral closures by \cite[Theorem 1.5]{Bha12} and then it descends via a $+$-splitting, compare with the discussion around \cite[Lemma 6.9]{BS17}. For the canonical sheaf, we know that pushforward respects dualizing complexes by Grothendieck--Serre duality and the preceding higher vanishing of the structure sheaf. Then, it suffices to observe that $S$ is globally $+$-regular, and hence Cohen--Macaulay (alternatively, we could have invoked Grauert--Riemenschneider vanishing).
\end{proof}

\subsection{Perfect BSDH varieties}

In this subsection, we study a perfect analog of BSDH varieties. This is relevant for us, because the Witt flag variety lives naturally in the world of perfect (ind-)schemes.

\begin{definition}
	A perfect BSDH variety over a perfect field $k$ is a perfect variety $X$ endowed with an inductive structure depending on its dimension as follows. For $\dim(X) = 0$, we require $X = \operatorname{Spec}(k)$. For $\dim(X) > 0$, this structure is given by a Zariski locally trivial $\mathbb{P}^{1,\mathrm{pf}}_k$-fibration $p \colon X \to Y$ equipped with a distinguished section, such that the base $Y$ is also a perfect BSDH variety.
\end{definition}

For such an object of dimension $n$, we can injectively construct closed immersions $\iota_J \colon X_J \to X$ of perfect schemes for any subset $J \subset [n]:=\{1,\dots,n\}$ and perfectly smooth projections $p_m\colon X\to X_m:=X_{[m]}$ for all $m\leq n$.
Again, we get curves $C_i:=X_{\{i\}}\subset X$ that carry isomorphisms $C_i \simeq \mathbb{P}^{1,\mathrm{pf}}_k$ (unique up to $\varphi$-twisting). We can also compute the Picard group to be $\mathbb{Z}[p^{-1}]^{\oplus n}$ given by restriction to $C_i$ as in the previous subsection. The next result explains the relationship between BSDH $k$-varieties and their perfect counterparts.

\begin{proposition} \label{prop:depBSDH} For a perfect BSDH $k$-variety $X$, 
 there exists a BSDH $k$-variety $X_0$ such that $X_0^{\mathrm{pf}}=X$.
\end{proposition}

\begin{proof}
	This is similar to \cite[Proposition 3.4]{HZ20}. We prove it by induction on $\mathrm{dim}(X)$. Consider a BSDH $k$-variety $Y_0$ deperfecting $Y$. Choose an affine open cover $\{ U_i\}$ of $Y$ such that $X_i:=X\times_Y U_i$ is a split $\mathbb{P}^1_k$-fibration with a section $U_i \to X_i$. We fix an isomorphism $X_i \simeq U_i \times \mathbb{P}^{1,\mathrm{pf}}_k$ intertwinning the given section with the origin section. This gives rise to a \v{C}ech cocycle for our affine open cover with values in a Borel subgroup of $\mathrm{PGL}_{2,k}$. Twisting by a sufficiently large power of Frobenius $\varphi$, we may assume that the cocycle has values in the affine open cover $\{U_{i,0}\}$ of the deperfection $Y_0$. This gives rise to a $\mathbb{P}^1_k$-fibration $X_0\to Y_0$ equipped with a section $Y_0\to X_0$ and perfecting to $X\to Y$ and its section.
\end{proof}

We expect that there exist BSDH $k$-varieties that are not $\varphi$-split, let alone
globally $+$-regular. We are also unaware of a notion of global $+$-regularity within the perfect setup, which relates to a good extent to the lack of an invertible line bundle as a dualizing sheaf. In the absence of a global property to descend along Stein factorizations, we ask the following.

\begin{question}
    Let $X$ be a perfect BSDH $k$-variety, $\calL$ be a semi-ample line bundle on $X$ and $S$ be the Stein factorization of $\calL$. Is $S$ a $\varphi$-rational perfect $k$-variety?
\end{question}

If $X = S$, this is true since $X$ is perfectly smooth.
As a first step toward rationality of $X\to S$, we prove that the Stein factorization map $X\to S$ factors as a sequence of maps with very simple geometric fibers.

\begin{lemma}\label{lem_sequence_maps_p1_or_singletons}
     Let $X$ be a perfect BSDH $k$-variety, $\calL$ be a semi-ample line bundle on $X$ and $S$ be the Stein factorization of $\calL$. Then, there is a sequence $Y_n:=X\to Y_{n-1}\to \dots \to Y_0:=S$ of maps of projective perfect $k$-varieties such that the geometric fibers of $Y_i\to Y_{i-1}$ are singletons or $\bbP^{1,\mathrm{pf}}_k$.
\end{lemma}

\begin{proof}
    Let $\calA_i$ be the pullback along $X\to X_i$ of an ample line bundle. We set $\calM_i:=\calL\otimes \calA_1\otimes \dots \otimes \calA_i$ and define $Y_i$ to be the Stein factorization of $\calM_i$. We claim that the natural map $X\to X_i\times S$ factors through a closed immersion $Y_i\hookrightarrow X_i\times S$. In order to check the factorization, suppose that $C$ is a connected perfect $k$-curve on which $\calM_i$ is torsion, i.e.~its degree vanishes. We have an expression
    \begin{equation}
        \deg(\calM_i|_C)=
        \deg(\calL|_C)+
        \sum_{j=1}^i\deg(\calA_j|_C)
    \end{equation}
    with each term nonnegative by semi-ampleness, so each term vanishes. This exactly means that $C$ is contracted under the maps to $X_j$ for $j\leq i$ and $S$, as desired. This gives us a tower of maps $Y_n\to \dots \to Y_i\to \dots \to Y_0$, but we still need to compute their geometric fibers.

    Next, we check that the resulting map $Y_i\to X_i\times S$ is injective. This follows by reversing the previous reasoning: given a connected perfect $k$-curve mapping to $X_i$ and $S$ via a constant map, we deduce that it also maps constantly to $X_j$ for all $j\leq i$. This now means that $\calL$ and $\calA_j$ for $j\leq i$ have trivial degrees on $C$, so the same is true for $\calM_i$. It follows that the entire curve $C$ maps via a constant map to $Y_i$, yielding injectivity. Now, we may compute the geometric fibers of $Y_i\to Y_{i-1}$, which are connected, otherwise $X\to Y_{i-1}$ would have disconnected geometric fibers.
    On the other hand, they embed as subschemes of the geometric fibers of $X_i\times S\to X_{i-1}\times S$, and the latter coincide with $\bbP^{1,\mathrm{pf}}_k$ by construction. Thus, the only possibilities for the former are either a point or $\bbP^{1,\mathrm{pf}}_k$.
\end{proof}

\begin{corollary} \label{lem:higherVanishing}
Let $X$ be a perfect BSDH $k$-variety, $\calL$ a semi-ample line bundle on $X$ and $f\colon X\to S$ its Stein factorization. Then, the direct images $R^if_\ast\calO_X$ vanish for $i>0$. Moreover, the cohomology groups $H^i(S,\calO_S)$ vanish for $i>0$.
\end{corollary}

\begin{proof}
This is a classical argument and could be compared as well to the statement in \cite[Lemma 2.8]{CX25} for $\bbF_p$-sheaves. Consider the factorization of $f$ into the maps $g_j\colon Y_j\to Y_{j-1}$. The geometric fibers of $g_j$ are either singletons or $\bbP^{1,\mathrm{pf}}$, both of whose structure sheaves have trivial higher cohomology. By base change for perfect $k$-schemes, see \cite[Lemma 3.18]{BS17}, we deduce that 
$R^ig_{j,\ast}\calO_X$ vanishes for $i>0$. By composition, this shows that $R^if_\ast\calO_X$ vanishes for $i>0$. In particular, it follows that $H^i(S,\calO_S)=H^i(X,\calO_X)$ by degeneration of the Grothendieck--Leray spectral sequence. To show its vanishing, we may now use the tower of Zariski local $\mathbb{P}^{1,\mathrm{pf}}_k$-fibrations $X_j\to X_{j-1}$.
\end{proof}

\begin{remark}
If Grauert--Riemenschneider vanishing held in perfect geometry like it holds for finite type $\varphi$-split varieties, see \cite[Theorem 1.3.14]{BK07}, then Cohen--Macaulayness of $S$ would be a consequence of Corollary \ref{lem:higherVanishing}. In the cone case, this can be shown, see \cite[Proposition 5.18]{BBL+23}, but otherwise counterexamples have been found in \cite[Theorem 5.2]{BBK23}.
\end{remark}

\section{Affine flag varieties}

In this section, we specialize our discussion to Demazure varieties that appear in the context of $p$-adic loop groups. Throughout this section $k$ denotes an algebraically closed field of characteristic $p > 0$. Let $F$ be a complete discretely valued field with ring of integers $\calO$ and residue field $k$. Fix a connected reductive group $G$ over $F$ and a parahoric $O$-model $\calG$ in the sense of Bruhat--Tits \cite{BT84}.

\subsection{Affine Schubert varieties}
We introduce the affine Schubert schemes following \cite[\S 3.2]{AGLR22}, with some simplifications since we assume $k$ is algebraically closed.

Let $\Alg$ denote the category of perfect $k$-algebras. For $R \in \Alg$ let $W(R)$ be the ring of $p$-typical Witt vectors over $R$. The ring of $\calO$-Witt vectors over $R$ is defined as
\begin{equation}W_{\calO}(R) = \begin{cases} W(R) \otimes_{W(k)} \calO, & \text{char}(F) = 0 \\ R \: \widehat{\otimes}_k \:\calO, & \text{char}(F) = p. \end{cases}\end{equation}
Note that if $\text{char}(F) = p$ and $t \in \calO$ is a uniformizer, then $\calO \cong k[\![t]\!]$ and $W_{\calO}(R) \cong R[\![t]\!]$. 

We define the following two functors $\Alg \to \text{Grp}$,
\begin{equation}\mathsf{G}(R) = G(W_\calO(R) \otimes_{\calO} F), \quad \mathsf{P}(R) = \calG(W_{\calO}(R)).\end{equation}
The affine flag variety for $\calG$ is the \'etale quotient $\mathsf{G}/\mathsf{P}$.
Its functor is representable by an increasing union of perfections of projective $k$-schemes. Indeed, if $\text{char}(F) = p$, then $\mathsf{G}/\mathsf{P}$ is the perfection of the affine flag variety in the sense of \cite{PR08} (which admits a natural moduli problem for all $k$-algebras), and if $\text{char}(F) = 0$, we obtain the affine flag variety in the sense of \cite{Zhu17} whose representability was proved in \cite[Corollary 9.6]{BS17}.

The Schubert varieties for the parahoric group scheme $\calG$ arise as the $\mathsf{P}$-orbit closures inside $\mathsf{G}/\mathsf{P}$. 
As we explain now, these are enumerated via double cosets of the Iwahori--Weyl group and this combinatorics captures their closure relations. Let $\mathbf{f}$ be the unique facet in the Bruhat--Tits building $\mathscr{B}(G,F)$ whose connected stabilizer is $\calG(O)=\mathsf{P}(k)$.
Let $A \subset G$ be a maximal $F$-split torus whose apartment contains $\mathbf{f}$. 
The centralizer $T=Z_G(A)$ is a maximal $F$-torus and we let $\calT$ be its connected N\'eron $O$-model. The Iwahori--Weyl group associated with $A$ is $W_{\mathsf{G}}: = N(F)/ \calT(\mathcal{O})$. The choice of an alcove $\mathbf{a}$ in the apartment of $A$ gives rise to a split exact sequence 
\begin{equation}1 \to W_{\aff} \to W_{\mathsf{G}} \to \pi_1(G)_I \to 1\end{equation} where $\pi_1(G)$ is the algebraic fundamental group and $I$ is the inertia group of $F$.
The affine Weyl group $W_{\aff}$ is the Coxeter group generated by the reflections in the walls of $\mathbf{a}$. By declaring elements of $\pi_1(G)_I$ to have length zero, $W_{\mathsf{G}}$ is a quasi-Coxeter group.  We denote the length function by $\ell$.
Let $W_{\mathsf{P}} \subset W_{\mathsf{G}}$ be the subgroup generated by reflections in the walls of $\mathbf{f}$. Then we have the Bruhat decomposition
\begin{equation}
	\mathsf{P}(k) \backslash \mathsf{G}(k) /\mathsf{P}(k) = W_{\mathsf{P}}\backslash W_{\mathsf{G}}/W_{\mathsf{P}}
\end{equation}
describing the $k$-valued points of the Hecke stack $\mathrm{Hk}_{\mathcal{G}}:=\mathsf{P}\backslash \mathsf{G}/\mathsf{P}$. Since these capture the entirety of the $\mathsf{P}$-orbits, we can now give the formal definition of Schubert varieties.

\begin{definition}
	Let $w \in W_{\mathsf{P}}\backslash W_{\mathsf{G}}/W_{\mathsf{P}}$. The affine Schubert variety $\mathsf{S}_w \subset \mathsf{G}/\mathsf{P}$ is the closure of the $\mathsf
	P$-orbit of any choice of lift of $w$ to $\mathsf{G}(k)/\mathsf{P}(k)$.
\end{definition}

The Schubert variety $\mathsf{S}_w$ of the affine group $\mathsf{G}$ is isomorphic to the perfection of a projective $k$-scheme. It has finitely many $\mathsf{P}$-orbits of the form $\mathsf{P}\cdot v$ for $v\leq w$ in the Bruhat order $\leq$. There is a refinement of this collection of closed subschemes obtained as $\mathsf{B}$-orbit closures in $\mathsf{G}/\mathsf{P}$, where $\mathsf{B}$ is the positive loop group of an Iwahori model $\mathcal{I}$ dilated from $\mathcal{G}$: these perfect $k$-varieties are in bijection with $ W_{\mathsf{G}}/W_{\mathsf{P}}$.

There is a general notion of convolution in this setup that will lead to a perfect BSDH $k$-variety. Let $w_\bullet=(w_1,\ldots,w_n)$ be a sequence of elements in $W_{\mathsf{G}}$. We define the convolution Schubert variety
\begin{equation}
	\mathsf{S}_{w_\bullet}:=\mathsf{G}_{w_1} \times^{\mathsf{P}} \cdots \times^{\mathsf{P}}\mathsf{S}_{w_n},
\end{equation}
where $\mathsf{G}_w\subset \mathsf{G}$ is the pullback of the Schubert variety $\mathsf{S}_w\subset \mathsf{G}/\mathsf{P}$ along the natural projection $\mathsf{G} \to \mathsf{G}/\mathsf{P}$ and the notation $\times^{\mathsf{P}}$ stands for the étale quotient by the diagonal $\mathsf{P}$-action on the adjacent factors. Any convolution Schubert variety admits a proper birational cover given by a Demazure variety that we introduce next, under the natural multiplication map. 

 If $\mathcal{G}=\mathcal{I}$ is an Iwahori and all the $w_i=:s_i$ have length $1$, then $\mathsf{G}_{s_i}$ identifies with the jet group $\mathsf{P}_{s_i}$ of the unique parahoric $O$-model $\mathcal{G}_{s_i}$ of $G$ such that $\mathcal{G}_{s_i}(O)$ stabilizes the codimension $1$ subfacet $\mathbf{f}_i \subset \bar{\mathbf{a}}$ fixed under $s_i$. Thus, we can write
\begin{equation}
	\mathsf{S}_{s_\bullet} = \mathsf{P}_{s_1} \times^{\mathsf{B}} \cdots \times^{\mathsf{B}} \mathsf{P}_{s_n}/\mathsf{B}.
\end{equation}
and call this a Demazure variety. Note that this is a perfect BSDH $k$-variety in our general sense. Indeed, $\mathsf{S}_{s_i}=\mathsf{P}_{s_i}/\mathsf{B}\simeq \mathbb{P}^{1,\mathrm{pf}}_k$ and there is a natural map $\mathsf{S}_{s_\bullet}\to \mathsf{S}_{t_\bullet}$ with $t_\bullet:=s_{\bullet <n}$ that is a locally trivial $\mathbb{P}^{1,\mathrm{pf}}_k$-fibration with a section induced by $\mathsf{B}\subset \mathsf{P}_{s_n}$. In particular, by Proposition \ref{prop:depBSDH} it is a perfectly smooth variety of dimension $n$ having a BSDH $k$-variety as a deperfection and its Picard group identifies with $\mathbb{Z}[p^{-1}]^{\oplus n}$ via the degree map, compare also with \cite[Lemma 4.8]{FHLR25} and \cite[Theorem 3.8]{AGLR22}. Usually there is some ambiguity in the perfect degree map that is resolved by fixing the deperfection. In this setup, however, we have a canonical deperfection of the perfect $k$-curve $\mathsf{S}_i$ given by the special fiber of the Iwahori $\mathcal{I}_k$, see \cite[Definition 3.14]{AGLR22}, so we fix that normalization. By \cite[Theorem 3.1]{HZ20}, the line bundle $\mathcal{O}(1,\dots,1)$ on $\mathsf{S}_{s_\bullet}$ is semi-ample, as it comes from an ample line bundle on $\mathsf{G}/\mathsf{B}$ via pullback. Later on, we will also need the following more general result.

\begin{lemma} \label{lem:amples}
	A line bundle $\calL$ on the Demazure variety $\mathsf{S}_{s_\bullet}$ is ample (resp.~semi-ample) if and only if $\mathrm{deg}(\calL)$ is a sequence of positive (resp.~non-negative) rationals and the subsequence indexed by any $s$ in $s_\bullet$ is strictly decreasing (resp.~decreasing).
\end{lemma}

\begin{proof}
	This essentially follows from \cite[Theorem 3.1]{HZ20}, but we give a self-contained proof. For the forward direction, notice that ampleness (resp.~semi-ampleness) is preserved under pull-back along a closed immersion (resp.~an arbitrary map). By restricting to $\mathsf{S}_s$, it follows that $\mathrm{deg}(\calL)$ consists of positive (resp.~non-negative) rationals if $\calL$ is ample (resp.~semi-ample). In order to obtain the monotonicity condition, we restrict $\calL$ to the convolution $\mathsf{S}_{(s,s)} $, and identify it with the usual product $\mathsf{S}_s^2$ via the first projection and multiplication. Observe that this isomorphism maps $\mathsf{S}_{(s,1)} \subset \mathsf{S}_{(s,s)} $ (resp.~ $\mathsf{S}_{(1,s)} \subset \mathsf{S}_{(s,s)} $) to the diagonal (resp.~second factor) of the untwisted product $\mathsf{S}_{s}^2$, so the claim is clear.
	
	For the converse, we consider the natural embedding $\mathsf{S}_{s_\bullet} \subset(\mathsf{G}/\mathsf{B})^n$ whose $i$-th coordinate is the multiplication of the first $i$ factors. Notice that $\calL$ is the pullback of a line bundle on the right side satisfying the same positivity condition on degrees. We may also pass to the adjoint quotient of $G$ and then to each of its simple $F$-factors, and hence assume that $G$ is an almost simple $F$-group. Now, we consider the closed embedding $\mathsf{G}/\mathsf{B}\to \prod_{s}\mathsf{G}/\mathsf{P}^s$, where $s$ runs through all simple reflections and $\mathsf{P}^s$ is the positive loop group in the unique maximal parahoric $O$-model such that $\calI(O)\subset \calG^s(O) $ but $s \notin \calG^s(O)$. Since $\mathrm{Pic}(\mathsf{G}/\mathsf{P}^s)=\bbZ[1/p]$, positivity equals ampleness for these partial flag varieties and the result is clear by pullback.
\end{proof}

Now, we can state our main result on global $\varphi$-regularity of Schubert varieties within a certain small range. It will be deduced from the corresponding global $\varphi$-regularity result for based BSDH $k$-varieties in the previous section, see \Cref{thm_globally_+_regular_abstract_demazure}.

\begin{theorem}
If the $\mathsf{P}$-action factors through the restriction of scalars of the mod $p$ fiber $\calG_{p=0}$ of the parahoric model, then $\mathsf{S}_{w_\bullet}$ has a globally $\varphi$-regular deperfection $\mathsf{S}_{w_\bullet}^{\mathrm{can}}$ that is compatibly $\varphi$-split with the corresponding deperfections for $v_\bullet \leq w_\bullet$.
\end{theorem} 

\begin{proof}
    Let $s_\bullet$ be a concatenation of reduced words for each $w_i$ such that the resulting map $\mathsf{S}_{s_\bullet}\to \mathsf{S}_{w_\bullet}$ is a proper birational cover, compare with \cite[Theorem 2.10]{CX25}.
    Recall that in \cite[Lemma 3.15]{AGLR22} it is explained how to construct a BSDH deperfection $\mathsf{S}^{\mathrm{can}}_{s_\bullet}$ by induction under the assumption on the $\mathsf{P}$-action: indeed, one checks that there are no $\varphi$-twists involved in the construction. This means that the associated line bundle $\mathcal{O}(1,\dots,1)$ of the BSDH $k$-variety $\mathsf{S}_{s_\bullet}^{\mathrm{can}}$ defined as in the previous section coincides after pulling back to $\mathsf{S}_{s_\bullet}$ with the semi-ample line bundle with degree sequence constant equal to $1$ mentioned right before stating the theorem.
    
    Therefore, we conclude that this deperfection $\mathsf{S}_{s_\bullet}$ is globally $\varphi$-regular by \Cref{thm_globally_+_regular_abstract_demazure} and \Cref{cor--splitBSDH}, and it is also compatibly $\varphi$-split with its BSDH subvarieties. Finally, we descend this to a deperfection of $\mathsf{S}_{w_\bullet}$  via proper pushforward as in \Cref{Fregular-prop_abs}. Note that the transition maps $\mathsf{S}_{v_\bullet}\to \mathsf{S}_{w_\bullet}$ are closed immersions since each twisted factor embeds in $\mathsf{G}/\mathsf{P}$.
\end{proof}

\begin{remark} 
	In equicharacteristic, the assumption on the $\mathsf{P}$-action is trivially satisfied and the $\mathsf{S}_{w_\bullet}^{\mathrm{can}}$ are the seminormalizations of the usual finite type Schubert $k$-varieties, so this gives a new and uniform proof of a result due to \cite[Theorem 1.4]{Cas22} for split $G$ and \cite[Theorem 4.1]{FHLR25} for general $G$. In mixed characteristic, the assumption holds for the $\mu$-admissible locus $\mathsf{A}_{G,\mu}$ defined in the next subsection. The known results in \cite[Lemma 3.15]{AGLR22} were derived from a comparison theorem with equicharacteristic and had an exception when $p=2$ and $G$ is odd unitary. In sum, we present a uniform proof of global $\varphi$-regularity in the expected range.
\end{remark}

We can also compute the Picard group of $\mathsf{S}_{w}^{\mathrm{can}}$ as follows: 

\begin{corollary}
	If the $\mathsf{P}$-action on $\mathsf{S}_w$ factors through the restriction of scalars of $\mathcal{G}_{p=0}$, then the natural map $\mathrm{Pic}(\mathsf{S}_{w}^{\mathrm{can}})\to \prod_{s\leq w}\mathrm{Pic}(\mathsf{S}_{s}^{\mathrm{can}})$ is an isomorphism.
\end{corollary}
\begin{proof}
	Recall by \Cref{lem_sequence_maps_p1_or_singletons} that $\mathsf{S}_{s_\bullet}\to \mathsf{S}_{w}$ factors as a sequence of maps whose non-trivial geometric fibers are isomorphic to $\bbP^{1,\mathrm{pf}}$. By compatible $\varphi$-splitness of the $\mathsf{S}_{v}^{\mathrm{can}}$, we can now write the deperfected Demazure resolution $\mathsf{S}_{s_\bullet}^{\mathrm{can}}\to \mathsf{S}_{w}^{\mathrm{can}}$ itself as a composition of rational maps with non-trivial fibers equal to $\mathbb{P}^1_k$, just like in \cite[Lemma 4.5]{FHLR25}. Then, one shows also inductively that a line bundle on $\mathsf{S}_{s_\bullet}^{\mathrm{can}}$ coming from the right side of the claimed equality pushes forward along the resolution to a line bundle on $\mathsf{S}_{w}^{\mathrm{can}}$, see \cite[Lemma 4.20]{FHLR25}.
\end{proof}

A consequence of this result is that we can now construct explicit $\varphi$-splittings of $\mathsf{S}_{s_{\bullet}}^{\mathrm{can}}$ via the Mehta--Ramanathan criterion in \cite{MR85}. We can even characterize their descents to the canonical deperfections of $\mathsf{S}_{w}\subset \mathsf{G}/\mathsf{B}$.

\begin{proposition}
	Assume $\mathsf{B}=\mathsf{P}$ and that the $\mathsf{B}$-action on $\mathsf{S}_w$ factors through the restriction of scalars of $\mathcal{I}_{p=0}$. Each line in $H^0(\mathsf{S}_{w}^{\mathrm{can}},\calO(1,\dots,1))$ spanned by a global section $\vartheta$ not vanishing at the origin determines a unique $\varphi$-splitting of $\mathsf{S}_{w}^{\mathrm{can}}$ compatibly splitting $\partial\mathsf{S}_{w}^{\mathrm{can}}$ and $\theta:=\mathrm{div}(\vartheta)$.
\end{proposition}
\begin{proof}
	Let us first show existence of $\vartheta$. We note every ample line bundle on $\mathsf{S}_{w}^{\mathrm{can}}$ is globally generated. Indeed, restriction to the origin yields a surjection upon taking global sections at the perfect level by the argument in \cite[Proposition 10.5]{BS17}. By the existence of some compatible $\varphi$-splitting, we get surjectivity as well at the level of canonical deperfections. Furthermore the base locus is $\mathsf{B}(k)$-equivariant and open, so it equals the whole space. 
	
	Let us now fix a reduced word $s_\bullet$ for $w$ and consider the pullback $f^\ast \vartheta$ along $f\colon \mathsf{S}^{\mathrm{can}}_{ s_{\bullet}}\to \mathsf{S}_{w}^{\mathrm{can}}$. We may now apply the Mehta--Ramanathan criterion, see \cite[Proposition 1.3.11]{BK07}, to produce a $\varphi$-splitting given by the $(p-1)$-th power of a global section $\sigma$ of the anticanonical sheaf such that $\mathrm{div}(\sigma)=\partial\mathsf{S}^{\mathrm{can}}_{s_{\bullet}}+\mathrm{div}(f^*\vartheta)$. We can descend this $\varphi$-splitting to $\mathsf{S}_{w}^{\mathrm{can}}$ and we immediately get that it compatibly splits $\partial\mathsf{S}_{s_\bullet}^{\mathrm{can}}$ and $\theta=\mathrm{div}(\vartheta)$ (which is therefore reduced).
	
	In order to prove the claim, we still have to characterize the descended $\varphi$-splitting on $\mathsf{S}_{w}^{\mathrm{can}}$ in unequivocal fashion, i.e., by proving that it is given by a specific $(p-1)$-th global section $\tau$ of the anticanonical sheaf of $ \mathsf{S}_{w}^{\mathrm{can}}$. 
	For this, we perform a calculation away from codimension $2$. Let $v \leq w$ be such that $\ell(v)=\ell(w)-1$. The Bruhat decomposition together with normality of the deperfections implies that the open locus of $\mathsf{S}_{w}^{\mathrm{can}}$ consisting of the $\mathsf{B}(k)$-orbits around $w$ and $v$ lifts isomorphically to $\mathsf{S}_{s_\bullet}^{\mathrm{can}}$. Now, one checks easily that the claimed $\tau$ exists and $\mathrm{div}(\tau)=\partial\mathsf{S}_{w}^{\mathrm{can}}+\mathrm{div}(\theta)$.
\end{proof}

Before moving on, we need to introduce the central extension $\hat{\mathsf{G}}$ of the loop group building on \cite[\S 4.1.3]{FHLR25}. This relates to the failure of $\mathsf{G}$-equivariance for line bundles on $\mathsf{G}/\mathsf{B}$. We assume from now on that $G$ is simply connected and almost simple. While one can show that every line bundle is $\mathsf{G}(k)$-equivariant, the Picard group of the stack $\mathsf{G}\backslash \mathsf{G}/\mathsf{B}$ identifies with the character group $X^*(\mathsf{T})$ of the maximal perfect split $k$-torus of $\mathsf{B}$ such that $\mathsf{T}(k)\subset \calT(O)$: this group has strictly lower rank than $\mathrm{Pic}(\mathsf{G}/\mathsf{B})$. 
Following \cite[Lemmas 4.10 and 4.12]{FHLR25}, the line bundle $\calO(\nu)$ associated with $\nu \in X^*(\mathsf{T})$ has degree $\langle a_s^\vee,\nu \rangle$ on $\mathsf{P}_s/\mathsf{B}$, where $s \in W_{\mathrm{af}}$ is a simple reflection and $a_s$ denotes the euclidean root underlying the associated affine root $\alpha_s$. The map $X^*(\mathsf{T})\to \mathrm{Pic}(\mathsf{G}/\mathsf{B})$ is a split injection of $\bbZ[p^{-1}]$-modules and its cokernel is a free $\mathbb{Z}[p^{-1}]$-module of rank $1$. Employing the same ordering as in \cite[Lemma 4.13]{FHLR25}, the line bundle $\calL:=\mathcal{O}(0,\dots,0,1)$ defines a section and pins down the central charge $c\colon \mathrm{Pic}( \mathsf{G}/\mathsf{B}) \to \mathbb{Z}[1/p]$ with kernel equal to $X^*(\mathsf{T})$. We define the central extension
\begin{equation}
	1 \to \bbG_{m,k}^{\mathrm{pf}} \to \widehat{\mathsf{G}} \to \mathsf{G} \to 1
\end{equation}
that classifies isomorphisms between $\calL$ and $g^*\calL$ for all $g\in \mathsf{G}(R)$ and every perfect $k$-algebra $R$. We could have chosen any other $\calL$ of central charge equal to $1$ and the resulting central extension would be independent of that choice by the proof of \cite[Lemma 4.27]{FHLR25}. We also denote by $\widehat{\mathsf{B}}$ the preimage of $\mathsf{B}$ along $\widehat{\mathsf{G}}\to \mathsf{G}$: we have $\widehat{\mathsf{G}}/\widehat{\mathsf{B}}=\mathsf{G}/\mathsf{B}$ and it follows that every line bundle thereon is now $\widehat{\mathsf{G}}$-equivariant. We employ similar notation for the preimage $\widehat{\mathsf{P}}$, resp.~$\widehat{\mathsf{T}}$, of $\mathsf{P}$, resp.~$\mathsf{T}$, along $\widehat{\mathsf{G}}\to \mathsf{G}$ and this does not change the corresponding flag variety. Using this, we define affine coroots.

\begin{definition}
	Let $s$ be a simple affine reflection in $W_{\text{af}}$ and $\calG_s$ be the associated parahoric group scheme. The affine coroot $\alpha_s^\vee \in X_*(\widehat{\mathsf{T}})$ is the coroot of the perfectly reductive quotient of $\widehat{\mathsf{P}}_s$ associated with the root $a_s \in X^*(S)$.

\end{definition}

If $\omega_s$ denotes the dual weights to the coroots $\alpha_s^\vee$, then $\calO(\omega_s)$ has degree $\delta_{ss'}$ on $\mathsf{P}_{s'}/\mathsf{B}$. If $\rho$ denotes the sum of the dual weights $\sum_s\omega_s$, then $\calO(\rho)=\mathcal{O}(1,\dots,1)$, compare with \cite[Lemma 4.17]{FHLR25}. Finally, we remark that there is a natural $\bbZ$-lattice
\begin{equation}
   \bbZ^{\mathrm{rk}G+1}\subset\bbZ[p^{-1}]^{\mathrm{rk}G+1} \simeq \mathrm{Pic}(\mathsf{G}/\mathsf{B})
\end{equation} 
spanned by the dual weights. We want to refer to these weights as integral and denote the previous inclusion by the shorthand $\Lambda_\bbZ \subset \Lambda_{\bbZ[p^{-1}]}$.

Now, we can formulate the Demazure character formula for general $G$. The central extension $\widehat{\mathsf{G}}$ acts on any line bundle $\calO(\nu)$. Assume $s_\bullet$ is a reduced word and the $\mathsf{B}$-action on $\mathsf{S}_{s_\bullet}$ factors through the restriction of scalars of $\calI_{p=0}$, so that we have a canonical deperfection $\mathsf{S}^{\mathrm{can}}_{s_\bullet}$. Moreover, assuming that $\nu\in\Lambda_\bbZ$,
then the line bundle $\calO(\nu)$ descends to $\mathsf{S}^{\mathrm{can}}_{s_\bullet}$. Hence, its cohomology groups have an associated character counting dimensions of affine integral weight spaces. This is an element of the group ring $\bbZ[\Lambda_\bbZ]=\oplus_{\nu\in\Lambda_\bbZ} \mathbb{Z}e^\nu$ of affine integral weights. We use exponential notation for this group ring to avoid confusion with sums of coefficients and sums of weights (which correspond to multiplication in the ring).

\begin{corollary}\label{cor_demazure_char_formula}
	Let $\nu \in\Lambda_\bbZ^+$ be a dominant integral weight. If the $\mathsf{B}$-action on $\mathsf{S}_{s_\bullet}$ factors through the restriction of scalars of $\mathcal{I}_{p=0}$, then
	\begin{equation}
		\mathrm{char}H^0(\mathsf{S}^{\mathrm{can}}_{s_\bullet},\calO(\nu))=\calD_{s_1}\circ\dots\circ \calD_{s_n}(e^\nu)
	\end{equation}
	where $\calD_s(e^\nu)=(1-e^{-a_s})^{-1}(e^\nu-e^{\nu-\langle \alpha_s^\vee,\nu+\rho\rangle a_s})$ is the Demazure operator.\end{corollary}

\begin{proof}
	Following the inductive computation in \cite[Theorem 7]{Lit98} identifies the right side with the character of the Euler characteristic $\chi(\mathsf{S}_{s_\bullet}^{\mathrm{can}}, \calO(\nu))$. 
The global $\varphi$-regularity of $\mathsf{S}^{\mathrm{can}}_{s_\bullet}$ and \cite[Theorem 1.6]{Bha12} applied to the semi-ample line bundle $\calO(\nu)$, compare with \Cref{lem:amples}, yields vanishing of higher cohomology, so the equality follows.
\end{proof}

Note that, given a Demazure resolution $\mathsf{S}_{s_\bullet}^{\mathrm{can}}\to \mathsf{S}_{w}^{\mathrm{can}}$, we can descend the previous calculation to the $\mathsf{S}_w^{\mathrm{can}}$, as the higher direct images of the structure sheaf vanish.
Finally, we briefly discuss what should happen when the $\mathsf{P}$-action on $\mathsf{S}_w$ does not factor through the mod $p$ fiber of $\calG$. The main issue that we face is that the deperfection constructed via \Cref{prop:depBSDH} will generally require us to twist by some partial Frobenius and we keep track of these choices via a sequence $q_\bullet$ of powers of $p$. The resulting theta divisor $\theta_{s_\bullet,q_\bullet}$ in the sense of \Cref{sec_abstract_BSDH} will no longer have degree $1$ on every $\mathsf{S}_s$, but will instead satisfy $\mathrm{deg}(\theta_{s_\bullet,q_\bullet})=q_\bullet$. In particular, this will almost never be effective if $q_\bullet$ is not constant (excepting products of different simple reflections). Nonetheless, we conjecture that there should be a statement as follows at the perfect level.

\begin{conjecture}
    The perfect $k$-variety $\mathsf{S}_w$ is $\varphi$-rational.
\end{conjecture}

Here is some modest evidence for the conjecture. For simplicity, we only address the Iwahori level case.  We were unable to make it work for Stein factorizations of abstract BSDH varieties, because it is not clear that the elements in the tower of partial resolutions of \Cref{lem_sequence_maps_p1_or_singletons} are smoothly equivalent to lower dimensional Stein factorizations. 

\begin{lemma}
    Assume $\mathsf{B}=\mathsf{P}$ and let $w \in W_\aff$. If $\mathsf{S}_v$ is $\varphi$-rational for every $v<w$, then $\mathsf{S}_w$ is Cohen--Macaulay.

\end{lemma}
\begin{proof}
   Let $s_\bullet$ be a word of simple reflections whose product is a reduced decomposition of $w$. Set $v$ to be the product of the reduced word $s_{\bullet>1}$ obtained by omitting the first letter.
   We have a partial resolution $m\colon \mathsf{S}_{s,v}\to \mathsf{S}_w$ whose non-singleton fibers are isomorphic to $\bbP^{1,\mathrm{pf}}$. We claim that $\mathsf{S}_{s,v}$ is also $\varphi$-rational. One way to prove this is to note that $\mathsf{S}_{s,v}$ and $\mathsf{S}_{s} \times \mathsf{S}_v$ have a common perfectly smooth cover, and then use ascent and descent for $\varphi$-nilpotence along smooth maps. For the latter, one can use \cite[Theorem B]{KMPS23} for descent, and Theorem \ref{prop-Fp-simple} along with the fact that intermediate extensions commute with smooth pullback as proved in \cite[Theorem 2.16]{Cas22} for ascent. It follows that $\bbF_p[\ell(w)]$ is a simple perverse sheaf on $\mathsf{S}_{s,v}$.

   Using vanishing of higher cohomology of $\bbP^{1,\mathrm{pf}}_k$,  the derived pushforward $Rm_\ast\bbF_p[\ell(w)]$ identifies with $\bbF_p[\ell(w)]$ on $\mathsf{S}_w$, compare also with \cite[Lemma 2.8]{CX25}. It is easy to see that $\mathbb{F}_p[\ell(w)]$ on $\mathsf{S}_w$ lies in perverse degrees $\leq 0$. To see that it lies in perverse degrees $\geq 0$, fix a point $x \in \mathsf{S}_w$. Let $\Spec(R)$ be a strict henselization of the local ring at $x$, let $\overline{x}$ be the closed point, and let $i \colon \overline{x} \rightarrow \mathsf{S}_w$ be the inclusion. We must show that $Ri^! \mathbb{F}_p[\ell(w)]$ lies in degrees $\geq - \dim\overline{\{x\}}$. Since $*$- and $!$-pullback agree for \'etale maps, we may replace $\mathsf{S}_w$ with $\Spec(R)$. If $h \colon Y \rightarrow \Spec(R) \times_{\mathsf{S}_w} \mathsf{S}_{s,v}$ is the fiber of $\overline{x} \rightarrow \Spec(R)$ over $m$, then by base change as in \cite[Theorem 8.4.9]{FuEtale} (this applies in our perfect setup by passage to a deperfection of $\mathsf{S}_w$), we are reduced to showing that $H^n(Y, Rh^! \mathbb{F}_p[\ell(w)]) = 0$ for $n < - \dim \overline{\{x\}}$. 

   We know that $Y$ is either $\overline{x}$ or $\mathbb{P}^1_{\overline{x}}$. In the former case, we conclude immediately by $\varphi$-nilpotence of $\mathsf{S}_{s, v}$ (even Cohen--Macaulayness suffices here). If $Y \cong \mathbb{P}^1_{\overline{x}}$, then there is some nonempty open subscheme $j \colon U \rightarrow  Y$ such that $j^*Rh^! \mathbb{F}_p[\ell(w)]$ lies in degrees $\geq - \dim \overline{\{x\}}$ for the standard t-structure, again by $\varphi$-nilpotence of $\mathsf{S}_{s, v}$ (since $\dim \overline{Y} = \dim \overline{\{x\}} + 1$). Likewise, if $Z$ is the complement of $U$ in $Y$ (a possibly empty disjoint union of points), then the $!$-pullback of $Rh^! \mathbb{F}_p[\ell(w)]$ to $Z$ lies in degrees $\geq - \dim \overline{\{x\}} $ for the standard t-structure, by Cohen--Macaulayness. Applying $m_*$ to the excision sequence for the $!$-pullback of $Rh^! \mathbb{F}_p[\ell(w)]$ with respect to $Y = U \cup Z$ then gives the necessary vanishing. 
   \end{proof}

\subsection{Local models}
In this section, we return to the initial assumption on $G$, which is no longer assumed to be simply connected almost simple, but a general connected reductive group over a local field $F$. We briefly need to refer to the theory of v-sheaves as in \cite{Sch17,SW20}, but the reader may treat this as a black box. Consider the Beilinson--Drinfeld Grassmannian $\Gr_{\calG}$ in the sense of \cite{SW20} defined over the v-sheaf $\Spd O$. Its generic fiber is isomorphic to the $B_{\mathrm{dR}}^+$-affine Grassmannian and its special fiber equals the v-sheaf attached to the affine flag variety $\mathsf{G}/\mathsf{P}$.
Let $\mu$ be a geometric conjugacy class of coweights with reflex field $E$. Then, the affine Grassmannian $\Gr_{G,E}$ base changed to $E$ contains a closed subsheaf $\Gr_{G,\leq \mu}$ arising as the $G(B^+_{\mathrm{dR}})$-orbit closure of $\mu(\xi)$. Following \cite{AGLR22}, we define the v-sheaf local model $\mathsf{M}^v_{\mu}$ as the v-sheaf closure of $\Gr_{G,\leq \mu}$ inside $\Gr_{\calG,O_E}$.

If $F$ has characteristic $p$ or $\mu$ is minuscule, there exists by \cite[Theorem 1.1]{AGLR22} and \cite[Corollary 1.4]{GL24} a unique flat normal projective $O_E$-scheme $\mathsf{M}_{\mu}$ with reduced special fiber whose associated v-sheaf equals $\mathsf{M}_{\mu}^v$. We call it the scheme-theoretic local model and, with the single exception of wild odd unitary groups (so only when $p=2$), it was shown in the corresponding statements of \cite{AGLR22,GL24} (relying on \cite[Theorem 1.2]{FHLR25}) that this scheme has $\varphi$-split special fiber and is Cohen--Macaulay. We can now use our global +-regularity result to remove this assumption from the computation of the special fiber in \cite{AGLR22,FHLR25,GL24}.

Before we state and prove it, we treat the deperfection of the $\mu$-admissible locus $\mathsf{A}_{\mu} \subset \mathsf{G}/\mathsf{P}$. Recall that the $\mu$-admissible set $\mathrm{Adm}_\mu$ of Kottwitz--Rapoport \cite{KR00} consists of all elements $w \in W_{\mathsf{G}}$ bounded by the translation $t_\lambda$ for some representative of $\mu$. Then, the $\mu$-admissible locus $\mathsf{A}_{\mu}$ is the union of all Schubert varieties $\mathsf{S}_{w}\subset \mathsf{G}/\mathsf{P}$ as $w$ runs over all double cosets with lifts in $\mathrm{Adm}_\mu$ (after fixing an Iwahori $O$-model $\mathcal{I}$ mapping to $\mathcal{G}$). Below, we show strong structure results on the canonical deperfection of $\mathsf{A}_{\mu}$ in the sense of \cite[Definition 3.14]{AGLR22}, generalizing \cite[Theorem 3.16]{AGLR22}.

\begin{proposition}
	Assume $F$ has characteristic $p$ or $\mu$ is minuscule. Then, $\mathsf{A}_{\mu}$ has a unique $\varphi$-split deperfection $\mathsf{A}^{\mathrm{can}}_{\mu}$ admitting the deperfections $\mathsf{S}^{\mathrm{can}}_{t_{\bar \lambda}}$ as compatibly $\varphi$-split closed subschemes for every representative $\bar\lambda$ of $\mu$.	
\end{proposition}

\begin{proof}
	Note that $t_{\bar \lambda}$ satisfies the $\mathsf{P}$-action assumption by the proof of \cite[Lemma 3.15]{AGLR22}, so the statement is reasonable. Consider more generally a finite subset $W$ of $W_{\mathsf{P}}\backslash W_{\mathsf{G}}/W_{\mathsf{P}}$, whose elements satisfy the $\mathsf{P}$-action assumption. Then, we claim that the Schubert scheme $\mathsf{S}_{W}$ in the sense of \cite[Definition 3.6]{AGLR22} obtained as the union of all Schubert varieties $\mathsf{S}_{w}$ with $w \in W$, is the finite colimit (which is neither filtered nor sifted) of these exact subvarieties along the various natural inclusions. We prove this observation by double induction on the dimension and the number of irreducible components (regarded as lexicographically ordered pairs). If $W$ has a maximal element, this is obvious as the colimit has a final term. Otherwise, we write it as a union $W_1 \cup W_2$ of incomparable sets, and verify that $\mathsf{S}_{W}$ equals the coproduct $\mathsf{S}_{W_1}\sqcup_{\mathsf{S}_{W_{12}}} \mathsf{S}_{W_2}$, where $W_{12}$ contains the maximal elements below $W_1$ and $W_2$. Indeed, the obvious map is a universal homeomorphism and one has a Mayer--Vietoris short exact sequence that computes structure sheaves.

	Now, assume that every element $w\in W$ has a reduced word $s_\bullet$ defined with respect to a fixed Iwahori $\calI$ mapping to $\calG$, such that $\mathsf{S}_{s_\bullet}^{\mathrm{can}}$ is a BSDH $k$-variety. We define the deperfection $\mathsf{S}_{W}^{\mathrm{can}}$ as the analogous colimit of the $\mathsf{S}_{w}^{\mathrm{can}}$. This exists at least as a functor, and we claim that it is representable by a scheme, that the transition maps $\mathsf{S}_{V}^{\mathrm{can}}\to \mathsf{S}_{W}^{\mathrm{can}}$ are closed immersions, and there is a compatible $\varphi$-splitting of them all. We perform again the same double induction argument. Assume first that $W$ has a maximal element: then existence is clear, it contains all the $\mathsf{S}_{v}^{\mathrm{can}}$ as closed subschemes, and these are all compatibly $\varphi$-split. As for $V\leq W$ without a maximal element, one can show as in the previous paragraph that the colimit is realized inside $\mathsf{S}_{W}^{\mathrm{can}}$: for the pushout step, one has to use the fact that compatibly split subschemes have reduced intersection. If $W$ has no maximal element, write it as a union $W_1 \cup W_2$ of incomparable sets. Now, existence of the colimit follows by the existence of pushouts along closed immersions, see \cite[Tag 0E25]{StaProj}. If we had a $\varphi$-splitting of $\mathsf{S}_{W}^{\mathrm{can}}$ compatible with all the other maps $\mathsf{S}^{\mathrm{can}}_{V}$, then the claim would follow.

	In order to do that, we assume first that $\calG=\calI$ is an Iwahori. Then, the line bundle $\calO(1,\dots,1)$ is defined on $\mathsf{S}_{W}^{\mathrm{can}}$ and it carries an origin-avoiding global section, again by the same double induction procedure. This shows that we get a compatible $\varphi$-splitting of $\mathsf{S}_{W}^{\mathrm{can}}$ compatible with all the smaller $V\leq W$. If $\calG$ is not an Iwahori, then we consider the birational universal homeomorphism $g \colon \mathsf{S}_{U}^{\mathrm{can}}\to \mathsf{S}_{W}^{\mathrm{can}}$, where $U$ consists of the obvious lifts defined by the reduced words. It suffices to compute the pushforward of the structure sheaf under $g_*$: this commutes with finite colimits at the derived level and then one gets the desired result at the abelian level by higher vanishing of direct images.
\end{proof}

The compatible $\varphi$-splitting of $\mathsf{A}_{\mu}^{\mathrm{can}}$ allows us to compute its Picard group in the same manner as in \cite[Proposition 4.23, Corollary 5.8]{FHLR25}, which is given by the obvious $\mathbb{Z}$-lattice in the Picard group of the corresponding connected component of $\mathsf{G}/\mathsf{P}$. In particular, we may speak unambiguously of the central charge $c_\calL$ for some line bundle $\calL$ on $\mathsf{A}_{\mu}^{\mathrm{can}}$ by $\mathsf{B}$-equivariantly translating to the neutral component, which is isomorphic to $\mathsf{G}_{\mathrm{sc}}/\mathsf{P}_{\mathrm{sc}}$. This forces us to regard the central charge $c_{\calL}$ as a tuple of integers obtained by splitting the simply connected cover of $G$ into simple factors.

Next, we compute the dimensions of global sections of an ample $\calL$ on $\mathsf{A}_{\mu}^{\mathrm{can}}$. We can also define a finite type $E$-scheme $\Gr_\mu^{\mathrm{can}}$ representing the generic fiber of $\mathsf{M}^v_\mu$: if $F$ has characteristic $p$, then this is again the seminormalization of the usual Schubert variety; if $F$ is $p$-adic and $\mu$ is minuscule, this is just the $E$-descent of the classical flag variety $G/P_\mu$ classifying parabolics of type $\mu$. The result below answers the coherence conjecture of Pappas--Rapoport \cite{PR08} in the $p$-adic minuscule setting and in the equicharacteristic setting, which is new in both settings if $p=2$ and $G$ is a wild odd unitary group.

\begin{theorem}\label{thm_coh_conj_local_model}
	Assume $F$ has characteristic $p$ or $\mu$ is minuscule. For any ample line bundle $\calL$ on $\mathsf{A}_{\mu}^{\mathrm{can}}$, there is an equality
	\begin{equation}
		\mathrm{dim}_k H^0(\mathsf{A}_{\mu}^{\mathrm{can}},\calL)=\mathrm{dim}_E H^0(\Gr_{\mu}^{\mathrm{can}},\mathcal{O}(c_{\calL}))
	\end{equation}
	of dimensions of global sections of line bundles, where $c_{\calL}$ is its central charge.
\end{theorem}

\begin{proof}
	The first order of business is to calculate the Euler characteristic of any line bundle $\calL$ on $\mathsf{A}_{\mu}^{\mathrm{can}}$. Below, we give a closed formula for the Euler characteristic on more general unions $\mathsf{S}_{W}^{\mathrm{can}}$ in terms of its irreducible Schubert subvarieties. 
	\begin{equation}
		\chi(\mathsf{S}_{W}^{\mathrm{can}},\calL)=\sum_{w_0} \mu_W(w_0)\chi(\mathsf{S}_{w_0}^{\mathrm{can}},\calL)
	\end{equation}
	where $\mu_W(w_0)$ is the difference between the number of even and odd chains $w_0<w_1<\dots<w_n \leq W$. We note that this formula is rather inefficient, but it would be difficult to write an optimal formula, as one would have to somehow capture the optimal intersection and union patterns. We regard it as an analog of the Möbius inversion formula for posets, but, since we could not find an appropriate reference, we give an explicit proof.
	
	Again, one has to perform the same double induction procedure. If $W$ has a maximal element $w$, then there is a bijection between $w$-avoiding chains and $w$-ending positive length chains. Since their lengths have different parity, the coefficient $\mu_W(w_0)$ vanishes for $w_0\neq w$, and only the term $\chi(\mathsf{S}_{w}^{\mathrm{can}})$ indexed by $w$ survives. Now write $W$ as the union of two sets $W_1$ and $W_2$ and let $W_{12}$ be the set of maximal elements below $W_1$ and $W_2$. Then, there is an equality
	\begin{equation}
		\chi(\mathsf{S}_{W}^{\mathrm{can}},\calL)+\chi(\mathsf{S}_{W_{12}}^{\mathrm{can}},\calL)=\chi(\mathsf{S}_{W_1}^{\mathrm{can}},\calL)+\chi(\mathsf{S}_{W_2}^{\mathrm{can}},\calL)
	\end{equation} 
	by Mayer--Vietoris, and we need to match the chain terms appearing in the Möbius inversion formula. If a chain ends in $W_{12}$, then it is summed twice on both sides. Otherwise, the term is counted only once on both sides, so we indeed get an equality. 
	
	Together with the Demazure character formula in \Cref{cor_demazure_char_formula} and higher vanishing for ample $\calL$, we get a purely combinatorial formula for the left side. We remark that, for triality factors, the dependence on $\mu$ rather than its image in the inertia coinvariants appears trickier because the splitting field can have either degree $3$ or $6$, but one checks that the $A_3$-average and the $S_3$-average of an arbitrary coweight coincide. This means that we can reduce the proof of the equality to tame $G$ in equicharacteristic (possibly even $0$), so the result follows from \cite[Theorem 3]{Zhu14}, see also \cite[Theorem 2.1]{GL24} for another proof.
\end{proof}

We can deduce from the above an identification of the special fiber of scheme-theoretic local models.
\begin{corollary} \label{Cor:SpecFiber}
	The special fiber of $\mathsf{M}_{\mu}$ equals $\mathsf{A}_{\mu}^{\mathrm{can}}$.
\end{corollary}

\begin{proof}
	If $p>2$ or $\Phi_G$ is reduced, this is \cite[Theorem 1.2]{AGLR22} and \cite[Theorem 1.2]{FHLR25}. Then, \Cref{thm_coh_conj_local_model} states that the Hilbert polynomials of $\mathsf{A}_{\mu}^{\mathrm{can}}$ and $\Gr_{\mu}^{\mathrm{can}}$ coincide. By \cite[Corollary 1.4]{GL24}, the special fiber of $\mathsf{M}_{\mu}$ is reduced with weak normalization equal to $\mathsf{A}_{\mu}^{\mathrm{can}}$ (essentially by definition). By flatness, the special fiber also shares the same Hilbert polynomial as $\Gr_{\mu}^{\mathrm{can}}$. Since the weak normalization injects on structure sheaves, the quotient module has trivial Hilbert polynomial, so it necessarily vanishes, and this yields the claim.
\end{proof}

\begin{remark}
    After our paper was first uploaded, He--Schremmer--Yu proved that $\mathsf{M}_{\mu}$ is Cohen--Macaulay in the last open case of wild odd unitary groups for $p=2$, building on our Corollary \ref{Cor:SpecFiber} together with their combinatorial methods, see \cite[Theorem B]{HSY26}.
\end{remark}

\bibliography{biblio.bib}
\bibliographystyle{alpha}
	
\end{document}